\documentclass[final]{siamonline171218}
\usepackage{amssymb}
\usepackage[T1]{fontenc}
\usepackage[normalem]{ulem}
\usepackage{graphicx,epstopdf} 
\usepackage[caption=false]{subfig} 

\usepackage{lipsum}
\usepackage{amsfonts}
\usepackage{graphicx}
\usepackage{epstopdf}
\usepackage{algorithmic}
\ifpdf
  \DeclareGraphicsExtensions{.eps,.pdf,.png,.jpg}
\else
  \DeclareGraphicsExtensions{.eps}
\fi

\numberwithin{theorem}{section}

\newcommand{\TheTitle}{From computation to comparison of tensor decompositions} 
\newcommand{\TheAuthors}{I. Domanov  and L. De Lathauwer}

\headers{\TheTitle}{\TheAuthors}

\title{{\TheTitle}\thanks{Submitted to the editors DATE.
\funding{This work was funded by (1) Research Council KU Leuven: C1 project c16/15/059-nD; (2) the Flemish Government under the ``Onderzoeksprogramma Artifici\"ele Intelligentie (AI) Vlaanderen'' programme; (3) F.W.O.:  project  G.0830.14N, G.0881.14N, 
	G.0F67.18N (EOS SeLMA); (4)  EU: The research leading to these results has received funding from the European Research Council under the European Union's Seventh Framework Programme (FP7/2007-2013) / ERC Advanced Grant: BIOTENSORS (no.  339804). This paper reflects only the authors' views and the Union is not liable for any use that may be made of the contained information.
}}}

\author{
	Ignat Domanov\thanks{
		Group Science, Engineering and Technology, KU Leuven - Kulak,
		E. Sabbelaan 53, 8500 Kortrijk, Belgium and
		Dept. of Electrical Engineering  ESAT/STADIUS KU Leuven,
		Kasteelpark Arenberg 10, bus 2446, B-3001 Leuven-Heverlee, Belgium
		(\email{ignat.domanov@kuleuven.be}, \email{lieven.delathauwer@kuleuven.be}).}
	\and
	Lieven De Lathauwer\footnotemark[2]
}

\usepackage{amsopn}

\ifpdf
\hypersetup{
  pdftitle={\TheTitle},
  pdfauthor={\TheAuthors}
}
\fi	

\newsiamremark{expl}{Example}
\newcommand{\upN}{{\hat N}}
\newcommand{\Bdiag}{\operatorname{Bdiag}}
\newcommand{\fF}{\mathbb F}
\newcommand{\unf}[2]{{\mathbf  #1}_{({#2}^c;#2)}}
\newcommand{\unfset}[2]{{#1}_{({#2}^c;#2)}}
\newcommand{\unfgen}[2]{{#1}_{({#2}^c;#2)}}
\newcommand{\unfone}[1]{{\mathbf  #1}_{(2,3;1)}}
\newcommand{\unftwo}[1]{{\mathbf  #1}_{(1,3;2)}}
\newcommand{\unfthree}[1]{{\mathbf  #1}_{(1,2;3)}}
\newcommand{\unfgenone}[1]{{#1}_{(2,3;1)}}
\newcommand{\colvec}[2]{
		ind^{I_{{#1}_1}\times\dots\times I_{{#1}_{#2}}}_{i_{{#1}_1},\dots,i_{{#1}_{#2}}}
	                    }

\newcommand{\nullsp}[1]{\operatorname{Null}\left( #1 \right) }

\newcommand{\mixmatrix}{\mathbf G}
\newcommand{\mixentry}{\mathbf g}

\newcommand{\moden}[1]{\text{mode-}$#1$
}
\newcommand{\Moden}[1]{\text{Mode-}$#1$
}

\begin{document}

\maketitle

\begin{abstract}
  Decompositions of higher-order tensors into  sums of simple terms are ubiquitous. We show that in order to verify that two tensors are  generated by  the same (possibly scaled) terms it is not necessary  to compute the individual decompositions. In general the explicit computation of such a decomposition may have high complexity and can be ill-conditioned. We now show that under some assumptions the verification can be reduced to a comparison of both the column and row spaces of the corresponding matrix representations of the tensors. We consider rank-1 terms as well as  low multilinear rank terms (also known as block terms) and show that the number of the terms and their multilinear rank can be inferred as well. The comparison relies only on numerical linear algebra and can be done in a numerically reliable way. We also  illustrate how our results can be applied to solve a multi-label classification problem that  appears in the context of blind source separation.
\end{abstract}


\section{Introduction}
Decompositions of tensors of order $N$ (i.e., $N$-way arrays of  real or complex numbers) into a sum of simple terms are ubiquitous. 
The most common simple term is a rank-$1$ tensor, i.e. a nonzero tensor whose columns
(resp. rows, fibers, etc.) are proportional. The corresponding decomposition into a minimal number of terms is known as
Canonical Polyadic Decomposition (CPD). 

It is well-known that for $N=2$, that is, in the matrix case, the decomposition in a minimal number of rank-$1$ terms  is not unique unless the matrix itself is rank-$1$:  indeed, any factorization  $\mathbf A={\mathbf X}^{(1)}{\mathbf X}^{(2)T}$ with full column rank factors ${\mathbf X}^{(1)}=[\mathbf x^{(1)}_1\ \dots\ {\mathbf x}^{(1)}_R]$ and ${\mathbf X}^{(2)}=[\mathbf x^{(2)}_1\ \dots\ {\mathbf x}^{(2)}_R]$ generates a valid decomposition $\mathbf A = \mathbf x^{(1)}_1\mathbf x^{(2)T}_1+\dots+\mathbf x^{(1)}_R\mathbf x^{(2)T}_R$, where $R$  is the rank of $\mathbf A$,  and this decomposition is not unique. 
 On the other hand, if  ${\mathbf X}^{(1)}$ and/or ${\mathbf X}^{(2)}$  are subject to constraints (e.g., triangularity or orthogonality), then the decomposition can be unique, but from an application point of view the imposed constraints can be unrealistic and the rank-$1$ terms not interpretable as meaningful ``data components''.
In contrast, for $N\geq 3$, that is, in the higher order tensor case, the unconstrained CPD  is easily  unique (see, for instance, \cite{Nick2018,Nick2017,PartI,PartII} and the references therein). 
Its uniqueness properties make  the  CPD  a fundamental tool 
for unique retrieval of data components, latent variable analysis, independent component analysis, etc., with countless applications in
chemometrics \cite{Bro2009}, telecommunication, array processing, machine learning, etc. \cite{LievenCichocki2013,ComoJ10,Kolda,TensRev2017}.

The higher order setting actually allows the recovery of terms that are more general than rank-$1$ terms. A MultiLinear (ML) rank-$(L_1,L_2,\dots)$ term  is a tensor whose columns (resp. rows, fibers, etc.) form a matrix of rank $L_1$ (resp. $L_2$, $L_3$, etc.).
Like CPD, a decomposition  into a sum of ML rank-$(L_1,L_2,\dots)$  terms (also known as block term decomposition)  is  unique under reasonably mild assumptions 
(see \cite{LDLBTDPartII,BTD2paper,BTD1paper} and the references therein),
so that it has found applications in  wireless communication \cite{LDL2008}, 
blind signal separation \cite{LievenLrLr1,Otto2016}, etc.

Tensor decompositions can be considered as  tools for {\it data analysis} that allow one to break a single (tensor) data set into small interpretable components.  It is known that, in general, the explicit computation 
		 of the CPD and the decomposition  into a sum of ML rank-$(L_1,L_2,\dots)$  terms  
		  may   have high complexity and can be ill-conditioned \cite{Nickill1,Nickill3,Nickill2}. In other words, the mildness of the uniqueness conditions comes with a numerical and a computational cost.
 
 In this paper we consider tensor decompositions from a fundamentally new perspective that is  closer to {\em pattern recognition}. Namely, we consider the following ``tensor similarity'' problem:
 \begin{itemize}
 	\item 
  {\em 
    How to verify that two  $I_1\times\dots\times I_N$ tensors are  generated by  the same (possibly scaled)  rank-$1$  terms?}
 \item 
 {\em
 	 	More generally, how to verify that two  $I_1\times\dots\times I_N$ tensors are  generated by  the same (possibly scaled) ML rank-$(L_1,L_2,\dots)$  terms?}
\end{itemize}
For brevity, our presentation will be in terms of the more general variant. The simpler (C)PD variant will follow as a special case (see, for instance, \cref{Thm: The_very_first_theorem}).

 An obvious approach would be to compute the decompositions of all tensors  and then to compare  them. 
This has two drawbacks. First, as mentioned above, the explicit computation of the decompositions may   have high complexity and can be ill-conditioned. Second, the approach may  fail if the tensors    are generated by the same (possibly scaled) terms in cases where 
the decompositions are not unique. 
   
  In this paper we will not compute the tensor decompositions. We will pursue a different approach, starting from the following trivial observation: 
   if 
 	\begin{equation}
 	 \text{a tensor }\mathcal B \text{ is a sum of (possibly scaled) terms from the decomposition of a tensor }
 	 \mathcal A,\\
 	 \label{eq:BinA}
 	\end{equation}
 	then
 	\begin{equation}
 	\operatorname{col}(\unfset{\mathbf B}{S})\subseteq \operatorname{col}(\unfset{\mathbf A}{S}) \text{ for all proper subsets } S \text{ of } \{1,\dots,N\},  
 	\label{eq:manyinclusions}   
 	\end{equation}
 	where $col(\cdot)$ denotes the column space of a matrix,  $S^c$ denotes the complement of the set $S$, and 
 	$\unfset{\mathbf A}{S}$ denotes the  $ (\prod\limits_{n\in S^c} I_n ) \times (\prod\limits_{n\in S} I_n)$ matrix representation of $\mathcal A$ (see \cref{subsection:redundancy} for  a formal definition of $\unfset{\mathbf A}{S}$). Actually we will explain that \cref{eq:manyinclusions} implies \cref{eq:BinA}. A clear advantage of the  approach based on the implication \cref{eq:manyinclusions}$\Rightarrow$\cref{eq:BinA} is that the conditions in \cref{eq:manyinclusions} 
 	rely only on numerical linear algebra and can be verified in a numerically reliable way. On the other hand, it is not known whether \cref{eq:BinA}
 can always be replaced by \cref{eq:manyinclusions}.
 	
 	Hence, the first contribution of this paper is to show that \cref{eq:manyinclusions} implies \cref{eq:BinA}. As a matter of fact, we will show that  \cref{eq:BinA}    follows from just $N$ conditions in \cref{eq:manyinclusions}, namely from
 	the conditions
 	\begin{equation}
 	\operatorname{col}(\unfset{\mathbf B}{n}) \subseteq \operatorname{col}(\unfset{\mathbf A}{n}), \quad n\in\{1,\dots, N\},
 	\label{eq:AincludesBnew}
 	\end{equation}
 	and  that   the $\frac{I_1\cdots I_N}{I_n}\times I_n$ matrices $\unfset{\mathbf A}{n}$ and $\unfset{\mathbf B}{n}$ in \cref{eq:AincludesBnew} can be used to compute the number of terms in the decompositions of $\mathcal A$ and $\mathcal B$ as well as their multilinear ranks.
 	We also consider a more general case where the inclusions in \cref{eq:AincludesBnew} are only known to hold  for some $n$ in $\{1,\dots, N\}$. 
  	
  	It is worth noting that the conditions
 	\begin{equation}
 	\operatorname{row}(\unfset{\mathbf B}{n}) \subseteq \operatorname{row}(\unfset{\mathbf A}{n}), \quad n\in\{1,\dots, N\},
 	\label{eq:AincludesBnewrow}
 	\end{equation} 
in which $row(\cdot)$ denotes the row space of a matrix, are more relaxed than the conditions in \cref{eq:AincludesBnew} (see Statement \ref{item:auxlemma0} of \cref{lemma:moved} below) and in general do not imply \cref{eq:BinA}. For instance, if $\frac{I_1\cdots I_N}{I_n}\geq I_n$, then  the conditions 	$\operatorname{row}(\unfset{\mathbf B}{n}) =\operatorname{row}(\unfset{\mathbf A}{n})$ $(=\fF^{I_n})$,  $n\in\{1,\dots, N\}$ hold for any
	generic tensors $\mathcal A$ and $\mathcal B$ (no matter whether they are generated by the same (possibly scaled) terms or not). 
  	
 	The second contribution of this paper is to show that the remaining $2^N-2-N$ conditions in \cref{eq:manyinclusions} are redundant, i.e., that the $N$ conditions in  	\cref{eq:AincludesBnew} imply all $2^N-2$ conditions in  \cref{eq:manyinclusions}.
(A fortiori, \cref{eq:BinA} follows from the $N$ conditions in \cref{eq:AincludesBnew}, as mentioned under the ``first contribution'' above.) 	
 	
 	Prior work on tensor similarity is limited to \cite{Frederik2018}. Both the present paper and \cite{Frederik2018} originated from the technical report \cite{LDLtenssim}.
 	 	The theoretical contributions of \cite{Frederik2018} related to the implication \cref{eq:AincludesBnew}$\Rightarrow$ \cref{eq:BinA}   rely on   prior knowledge on the decompositions of $\mathcal A$ and $\mathcal B$\footnote{Namely, the working assumption in  \cite{Frederik2018} is that both tensors $\mathcal A$ and $\mathcal B$ admit decompositions
 	of the same type (CPD, decomposition in ML rank-$(L,L,1)$ terms, decomposition in ML rank-$(L,L,\cdot)$ terms), that the decompositions include the same number of terms, and that in the latter two decomposition types the terms of $\mathcal A$ and $\mathcal B$ can be matched so that their ML ranks are equal.}  and can be summarized as follows: if $N=3$ and \cref{eq:AincludesBnew} holds with ``$\subseteq$'' replaced by ``$=$'', then  $\mathcal A$ and $\mathcal B$ are generated by the same (possibly scaled) terms. The results obtained in the current paper imply that  the prior knowledge on the decompositions is not needed.	Further, \cite{Frederik2018} presents applications 
 		 in the context of emitter 	movement detection and fluorescence data analysis.

The  paper is organized as follows.  In \cref{subsec:basic_definitions,subsec:problem_statement} we introduce
tensor related notations and formalize the problem statement, respectively. \Cref{sec:preliminary} contains preliminary results. In \cref{sec:aux}, for the convenience of the reader, we remind the 
primary decomposition theorem and the Jordan canonical form. \Cref{secmovedlemma} contains  an auxiliary result  about the simultaneous compression of tensors $\mathcal A$ and $\mathcal B$ for which the first $\upN$ inclusions  in \cref{eq:AincludesBnew} hold (\cref{lemma:moved}). The  main results are given in \cref{sec:Main}. In \cref{sec:Mainconnections}
we establish connections between the terms in the decompositions of  tensors $\mathcal A$ and $\mathcal B$ that satisfy the  conditions
in \cref{eq:AincludesBnew}  (\cref{thm3rdorder,thm:newverygeneral theorem}). In \cref{subsection:redundancy} 
we show that the $N$ conditions in \cref{eq:AincludesBnew} imply the $2^N-2$ conditions in \cref{eq:manyinclusions} (\cref{corr:redandancy}).
In \cref{sec4:illustration} we illustrate how our results can be applied to solve a multi-label classification problem that appears in the context of blind source separation.  
\section{ Basic definitions and problem statement}
\subsection{Basic  definitions} \label{subsec:basic_definitions}
\subsubsection*{Matrix representations} Let $1\leq n\leq N$.
{\em A \moden{n} matrix representation} of a tensor $\mathcal A\in\fF^{I_1\times \dots \times I_N}$ is a matrix  $\unf{A}{n}\in\fF^{\frac{I_1\cdots I_N}{I_n}\times I_n}$  
whose columns are  the vectorized \moden{n} slices of $\mathcal A$. Using Matlab colon notation,
the columns of $\unf{A}{n}$ are the vectorized $I_1\times\dots\times I_{n-1}\times 1\times I_{n+1}\times\dots\times I_N$ tensors 
$\mathcal A(:,\dots,:,1,:,\dots,:), \dots, \mathcal A(:,\dots,:,I_n,:,\dots,:)$.
 Formally,
\begin{equation}\label{eq:mode_n-unfolding}
\text{the } (1+\sum_{\substack{k=1\\ k\ne n}}^N (i_k-1)\prod_{\substack{l=1\\ l\ne n}}^{k-1}I_l,i_n)\text{th entry of }\unf{A}{n} = \text{ the }(i_1,\dots, i_N)\text{th entry of }\mathcal A. 
\end{equation}
\subsubsection*{\Moden{n} product}
If for some tensor $\mathcal D\in\fF^{I_1\times \dots I_{n-1}\times L_n\times  I_{n+1}\times I_N}$ and  matrix $\mathbf X^{(n)}\in\fF^{I_n\times L_n}$,
\begin{equation}
\label{eq:An_one}
\unf{A}{n} = \unf{D}{n}\mathbf X^{(n)T},
\end{equation}
i.e., if the \moden{n} fibers of $\mathcal A$ are obtained by multiplying the corresponding \moden{n} fibers of $\mathcal D$ by $\mathbf X^{(n)}$, 
then we say that
 $\mathcal A$ is {\em  the \moden{n} product} of a  $\mathcal D$ and  $\mathbf X^{(n)}$ and write
$
\mathcal A=\mathcal D\bullet_n \mathbf X^{(n)}
$. 
It can be easily verified that the remaining  $N-1$ matrix representations of $\mathcal A$ can be factorized as
\begin{equation}
\label{eq:Ak_remaining}
\unf{A}{k} = \left(\bigotimes_{l=1,l\ne k}^{n-1} \mathbf I_{I_l}
\otimes \mathbf X^{(n)} \otimes  
\bigotimes_{l=n+1,l\ne k}^{N} \mathbf I_{I_l}\right) \unf{D}{k},\qquad k\in\{1,\dots,N\}\setminus \{n\}.
\end{equation}
where $\mathbf I_{I_l}$ and  ``$\otimes$'' denote the $I_l\times I_l$ identity matrix and the Kronecker product, respectively.
\subsubsection*{Several products in the same mode or across modes}
It easily follows from \cref{eq:An_one} that 
for compatible matrix and tensor dimensions, 
$$
\left(\left(\left(\mathcal D\bullet_n\mathbf X_1^{(n)} \right) \bullet_n \mathbf X_2^{(n)}\right)\dots \bullet_n \mathbf X_k^{(n)}\right)=
\mathcal D\bullet_n \left(  \mathbf X_k^{(n)}\cdots \mathbf X_1^{(n)}\right).
$$
Let $\upN\leq N$ and
\begin{equation}
\label{eq:DandXXXdims}
\mathcal D\in\fF^{L_1\times\dots \times L_{\upN}\times I_{\upN+1}\times\dots\times I_{N}},\qquad \mathbf X^{(1)}\in\fF^{I_1\times L_1},\dots, \mathbf X^{(\upN)}\in\fF^{I_{\upN}\times L_{\upN}}.
\end{equation}
For products across different modes, we have
\begin{align}
\mathcal D\bullet_1\mathbf X^{(1)}\dots\bullet_{ N} \mathbf X^{( N)} :=&\left(\left(\left(\mathcal D\bullet_1\mathbf X^{(1)} \right) \bullet_2 \mathbf X^{(2)}\right)\dots \bullet_{ N} \mathbf X^{( N)}\right)= \label{eq:defmanyprods}\\
&\left(\left(\left(\mathcal D\bullet_{i_1}\mathbf X^{(i_1)} \right) \bullet_{i_2} \mathbf X^{(i_2)}\right)\dots \bullet_{i_{ N}} \mathbf X^{(i_{ N})}\right)\nonumber
\end{align}
for any permutation $i_1,\dots,i_{N}$ of $1,\dots,N$.
It
 follows from \cref{eq:An_one,eq:Ak_remaining,eq:defmanyprods}, that the matrix representations of $\mathcal A=\mathcal D\bullet_1\mathbf X^{(1)}\dots\bullet_{N} \mathbf X^{(N)}$  are given by
\begin{equation}
\label{eq:Anmany}
\unf{A}{n}= \left(\bigotimes_{k=1,k\ne n}^N\mathbf X^{(k)}\right)
\unf{D}{n}\mathbf X^{(n)T},\qquad n\in\{1,\dots,N\}.
\end{equation}
If $\mathcal A=\mathcal D\bullet_1\mathbf X^{(1)}\dots\bullet_{\upN} \mathbf X^{(\upN)}$ with $\upN < N$, then the identities in \cref{eq:Anmany} hold with
$\mathbf X^{(\upN+1)}=\mathbf I_{I_{\upN+1}},\dots,\mathbf X^{(N)}=\mathbf I_{I_N}$. That is,
\begin{align}
\unf{A}{n}&= \left(\bigotimes_{k=1,k\ne n}^{\upN}\mathbf X^{(k)} \otimes \bigotimes_{k=\upN+1}^N\mathbf I_{I_k}\right)
\unf{D}{n}\mathbf X^{(n)T},\qquad n\in\{1,\dots,\upN\},\label{eq:identupN1}\\
\unf{A}{n}&= \left(\bigotimes_{k=1}^{\upN}\mathbf X^{(k)} \otimes \bigotimes_{k=\upN+1,k\ne n}^N\mathbf I_{I_k}\right)
\unf{D}{n},\qquad n\in\{\upN+1,\dots,N\}.\label{eq:identupN2}
\end{align}
\subsubsection*{ML rank of a tensor}
By definition, 
$$
\mathcal A\ \text{ is ML rank-}(L_{1},\dots,L_{\upN },\cdot,\dots,\cdot)\overset{\text{def}}{\Longleftrightarrow} r_{\unf{A}{n}}=L_{n},\quad n\in\{1,\dots,\upN\},\quad 2\leq \upN\leq N,
$$
that is, $L_{n}$ is the dimension of the subspace spanned by the \moden{n} fibers of $\mathcal A$.
It can be shown that $\mathcal A$ is ML rank-$(L_{1},\dots,L_{\upN },\cdot,\dots,\cdot)$ if and only if it admits the factorization
$\mathcal A=\mathcal D\bullet_1\mathbf X^{(1)}\dots\bullet_{\upN} \mathbf X^{(\upN)}$ such that
$\mathcal D$, $\mathbf X^{(1)},\dots, \mathbf X^{(\upN)}$ have dimensions as in \cref{eq:DandXXXdims} and
$\mathbf X^{(1)},\dots, \mathbf X^{(\upN)}, \unf{D}{1},\dots,\unf{D}{\upN}$ have full column rank. In this paper we assume that the tensor dimensions  have been permuted so that  we can  just specify the rank values for the first $\upN$ matrix representations  of $\mathcal A$.   A special case of the factorization
$\mathcal A=\mathcal D\bullet_1\mathbf X^{(1)}\dots\bullet_{\upN} \mathbf X^{(\upN)}$, where
$\upN =N$,  $\mathbf X^{(n)}$ equals the ``$U$'' factor in the compact Singular Values Decomposition (SVD) of $\unf{A}{n}$, 
and  $\mathcal D=\mathcal A\bullet_1\mathbf X^{(1)H}\dots \bullet_N\mathbf X^{(N)H}$ is known as the MLSVD  of $\mathcal A$  and is used for the compression of an $I_1\times\dots\times I_N$ tensor to the size $L_1\times\dots\times L_N$ \cite{2000LDLetall}. By setting  $\mathbf X^{(n)}$ equal to the identity matrix for $n=\upN+1,\dots,N$, we   compress only along  the first $\upN$ dimensions.
\subsubsection*{ML rank-$(L_{1r},\dots,L_{\upN r},\cdot,\dots,\cdot)$ decomposition of a tensor}
In this paper we consider the  decomposition of  $\mathcal A$ into a sum of ML rank-$(L_{1r},\dots,L_{\upN r},\cdot,\dots,\cdot)$ terms:
\begin{gather}
\mathcal A =\sum\limits_{r=1}^R  \mathcal D_r\bullet_1\mathbf X^{(1)}_r\dots\bullet_\upN\mathbf X^{(\upN)}_r,\quad 2\leq \upN\leq N,\label{eq:MLrankdecomp}\\
\mathcal D_r\in\fF^{L_{1r}\times\dots\times  L_{\upN_r}\times I_{\upN_r+1}\times\dots\times I_N}, \quad \mathbf X^{(n)}_r\in\fF^{I_n\times L_{nr}}, \quad n\in\{1,\dots, \upN\},\quad r\in\{1,\dots,R\}.\nonumber
\end{gather}
In our derivation we will also use a matricized version of \cref{eq:MLrankdecomp}. It can be obtained as follows.
First, we call
\begin{equation}
\label{eq:defXn1}
\mathbf X^{(n)}:=[\mathbf X^{(n)}_1\ \dots\ \mathbf X^{(n)}_R] \in\mathbb F^{I_n\times \sum\limits_{r=1}^R L_{nr}}, \qquad n\in\{1,\dots,\upN\},
\end{equation}
the {\em concatenated factor matrices of} $\mathcal A$.
If further we set
\begin{equation}
\label{eq:defXn2}
\mathbf X^{(n)}:=[\mathbf I_{I_n}\ \dots\ \mathbf I_{I_n}]\in\mathbb F^{I_n\times RI_n},\qquad n\in\{\upN+1,\dots,N\},
\end{equation}
then, by \cref{eq:Anmany}, we can express \cref{eq:MLrankdecomp} in a matricized way as
\begin{multline}
\label{eq:LL1matr}
\unf{A}{n}= \sum\limits_{r=1}^R   
\left(\bigotimes_{l=1,l\ne n}^N\mathbf X_r^{(l)}\right)
\unfgen{{\mathbf D}_r}{n}\mathbf X_r^{(n)T}=\\
\left(\bigodot\limits_{l=1,l\ne n}^N \mathbf X^{(l)}\right)\Bdiag(
\unfgen{{\mathbf D}_1}{n}, \dots, \unfgen{{\mathbf D}_R}{n})
\mathbf X^{(n)T},\qquad n\in\{1,\dots,N\},
\end{multline}
where 
\begin{equation}
\label{eq:bigodotXk}
\bigodot\limits_{l=1,l\ne n}^N \mathbf X^{(l)} :=\left[
\bigotimes_{l=1,l\ne n}^N\mathbf X_1^{(l)}\ \dots\ \bigotimes_{l=1,l\ne n}^N\mathbf X_R^{(l)}  \right]
\end{equation}
and
$\Bdiag(\unfgen{{\mathbf D}_1}{n}, \dots, \unfgen{{\mathbf D}_R}{n})$ denotes a block-diagonal matrix with  the matrices
$\unfgen{{\mathbf D}_1}{n}$, $\dots$, $\unfgen{{\mathbf D}_R}{n}$ on the diagonal.

	Note that  \cref{eq:MLrankdecomp} captures several well-studied decompositions as special cases (see also the introduction).	 
	If $\upN=N$ and $L_{1r}=\dots=L_{N r}=1$ for all $r$, then  all terms in  \cref{eq:MLrankdecomp} are rank-$1$ tensors, so 
	\cref{eq:MLrankdecomp} reduces to a polyadic decomposition of $\mathcal A$.
	It can easily be verified that if $\upN=2$, $N=3$, and $L_{1r}=1$ for all $r$, then  
	the ML rank-$(1,L_{2r},\cdot)$ terms in \cref{eq:MLrankdecomp} are actually 
	ML rank-$(1,L_{2r},L_{2r})$ terms. Thus, \cref{eq:MLrankdecomp} reduces to the decomposition into a sum of 	ML rank-$(1,L_{2r},L_{2r})$ terms.
	Finally, if  $\upN=2$ and $N=3$, then \cref{eq:MLrankdecomp} is a tensor reformulation of the joint block diagonalization problem. Namely,
	  \cref{eq:MLrankdecomp} means that the frontal slices   of $\mathcal A$ can   simultaneously be factorized as 
	$$
	\mathbf A(:,:,i) = \mathbf X^{(1)} \Bdiag({\mathbf D}_1(:,:,i),\dots, {\mathbf D}_R(:,:,i))  \mathbf X^{(2)T},\quad i=1,\dots,I_3,  
	$$
	where ${\mathbf D}_r(:,:,i)\in\mathbb F^{L_{1r}\times L_{2r}}$. 
	\subsection{Problem statement}\label{subsec:problem_statement}
 Assume that  a tensor  $\mathcal B\in\fF^{I_1\times \dots \times I_N}$ consists of the same 
 ML rank-$(L_{1r},\dots,L_{\upN r},\cdot,\dots,\cdot)$ terms as $\mathcal A$, but possibly differently scaled:
\begin{equation}
\mathcal B =\sum\limits_{r=1}^R  \lambda_r\mathcal D_r\bullet_1\mathbf X^{(1)}_r\dots\bullet_\upN\mathbf X^{(\upN)}_r,\qquad \lambda_1\cdots\lambda_R\ne 0.\label{eq:decompB}
\end{equation}
Then by \cref{eq:LL1matr},
\begin{multline}
\label{eq:matrunfB}
\unf{B}{n}= \left(\bigodot\limits_{k=1,k\ne n}^N \mathbf X^{(k)}\right)
\Bdiag(\lambda_1\unfgen{{\mathbf D}_1}{n}, \dots, \lambda_R\unfgen{{\mathbf D}_R}{n})
\mathbf X^{(n)T}=\\
\left(\bigodot\limits_{k=1,k\ne n}^N \mathbf X^{(k)}\right)
\Bdiag(\unfgen{{\mathbf D}_1}{n}, \dots, \unfgen{{\mathbf D}_R}{n})
\Bdiag(\lambda_1\mathbf I_{L_{n1}},\dots,\lambda_R\mathbf I_{L_{nR}})\mathbf X^{(n)T}.
\end{multline}
Assume that $\upN\geq 2$ and that the matrices
\begin{equation}
\label{eq:Xnhavefcr}
\mathbf X^{(1)},\dots,\mathbf X^{(\upN)}\ \text{ have full column rank.}
\end{equation}
It can be easily shown\footnote{Indeed, the result holds since, by assumption \cref{eq:Xnhavefcr},  the  first $\upN-1$ factors $\mathbf X^{(l)}$  have full column rank and, by construction, the remaining factors  do  not have zero columns.}  that the matrices in \cref{eq:bigodotXk} have full column rank for all $n$.
Hence, by \cref{eq:LL1matr,eq:matrunfB}, 
the column spaces of the first $\upN$ matrix representations of $\mathcal A$ and $\mathcal B$ coincide:
 \begin{equation}
 \label{eq:colAn=colBn}
 \operatorname{col}(\unf{A}{n})=\operatorname{col}(\unf{B}{n}),\qquad n\in\{1,\dots,\upN\}.
 \end{equation}
  If we further limit\footnote{\Cref{lemma:moved} below implies that assumption \cref{eq:Xnhavefcr} can always be replaced by assumption \cref{eq:Xnnonsingular}. Computationally, this can be done by Multilinear Singular Value Decomposition (MLSVD) \cite{2000LDLetall, Tucker1964, Tucker1966}.} ourselves to the case where
the matrices 
\begin{equation}
\label{eq:Xnnonsingular}
\mathbf X^{(1)},\dots,\mathbf X^{(\upN)}\ \text{ are square and nonsingular,}
\end{equation} 
then, obviously,
\begin{equation}
\unf{B}{n}=\unf{A}{n}\mathbf M_n,\qquad n\in\{1,\dots,\upN\},\label{eq:An=BnMn}
\end{equation}
where
\begin{equation}
 \mathbf M_n =  \left(\mathbf X^{(n)T}\right)^{-1} \Bdiag(\lambda_1\mathbf I_{L_{n1}},\dots,\lambda_R\mathbf I_{L_{nR}})\mathbf X^{(n)T},\qquad n\in\{1,\dots,\upN\}.\label{eq:MnEVD}
\end{equation}
Thus, if \cref{eq:MLrankdecomp,eq:decompB,eq:Xnnonsingular} hold, then the column spaces of the first $\upN$ matrix representations of $\mathcal A$ and $\mathcal B$ coincide, the matrices $\mathbf M_n:={{\mathbf A}^\dagger_{(n^c;n)}}\unf{B}{n}$ have the same spectrum $\lambda_1,\dots,\lambda_R\in\fF$ and can be diagonalized, $n=1,\dots,\upN$. Moreover,
the concatenated factor matrices $\mathbf X^{(n)}$  and the sizes of blocks $L_{nr}$ (and hence the overall decompositions of $\mathcal A$ and $\mathcal B$)  can be recovered from the EVDs of $\mathbf M_1,\dots, \mathbf M_\upN$.

In this paper we consider the inverse problem: we assume that the column spaces of the first $\upN$ matrix representations of $\mathcal A$ and $\mathcal B$ coincide and we investigate how the ML rank decompositions $\mathcal A$ and $\mathcal B$ relate to each other.
In particular, we obtain the following result.
\begin{theorem}\label{Thm: The_very_first_theorem}
Let $\mathcal A, \mathcal B\in\mathbb C^{I_1\times \dots \times I_N}$ and $2\leq\upN\leq N$. Assume that \cref{eq:Xnnonsingular,eq:An=BnMn} hold and that
at least one of the matrices $\mathbf M_1,\dots,\mathbf M_{\upN}$ can be diagonalized.
Then the following statements hold.
\begin{enumerate}
	\item
The matrices $\mathbf M_1,\dots,\mathbf M_{\upN}$ have the same spectrum.
\item
All matrices $\mathbf M_1,\dots,\mathbf M_{\upN}$ can be diagonalized.
\item
Let the distinct eigenvalues of $\mathbf M_n$ be $\lambda_1,\dots,\lambda_R$ with respective multiplicities
$L_{n1},\dots,L_{nR}$ and let $\mathbf X_n\in\mathbb C^{I_n\times I_n}$ be a nonsingular matrix such that
\cref{eq:MnEVD} holds.
Then  $\mathcal A$ and $\mathcal B$ admit the ML rank-$(L_{1r},\dots,L_{\upN r},\cdot,\dots,\cdot)$ decompositions in \cref{eq:MLrankdecomp,eq:decompB}, respectively. In particular, if $L_{nr}=1$ for all $n$ and $r$, then  
$\mathcal A$ and $\mathcal B$ are  generated by  the same  (possibly scaled) $R$  rank-$1$  terms.
\end{enumerate}
\end{theorem}
\begin{proof}
		The proof follows from \cref{thm:newverygeneral theorem} below.
	\end{proof}
The theorem can be used as follows. First, the matrices $\mathbf M_1,\dots,\mathbf M_{\upN}$ are found from the sets of linear equations \cref{eq:An=BnMn}.
(If any  of the sets of linear equations does not have a solution, then $\mathcal B$ is not of the form \cref{eq:decompB}, i.e., it cannot be generated by terms from the decomposition of $\mathcal A$.) The number of terms $R$ is found as the number of distinct eigenvalues of $\mathbf M_n$, $1\leq n\leq \upN$. The distinct eigenvalues themselves correspond to the scaling factors $\lambda_r$ in \cref{eq:decompB}. Both $R$ and the eigenvalues $\lambda_r$ are necessarily the same for all $\mathbf M_n$, but the multiplicities can be different. The multiplicity of $\lambda_r$ in the EVD of $\mathbf M_n$ corresponds to the $n$th entry $L_{nr}$ in the ML rank of the $r$th term. The larger $\upN$, the more the terms are specified. The minimal value for $\upN$ is $2$, since a decomposition in ML rank-$(L_{1r},\cdot,\dots,\cdot)$ terms is meaningless.

So far, we have explained the use of the theorem for decompositions that are exact. Obviously, the theorem also suggests a procedure for approximate decompositions (of noisy tensors). The equations in \cref{eq:An=BnMn} may be solved in least squares sense. The eigenvalues $\lambda_{nr}$ of the matrices
$\mathbf M_1,\dots,\mathbf M_{\upN}$ may be averaged over $n$ to obtain estimates of $\lambda_r$. The values $L_{nr}$, $1\leq n\leq \upN$, $1\leq r\leq R$ may be estimated by assessing how close the eigenvalues $\lambda_{rn}$ are to the averaged values $\lambda_r$.
\section{Preliminaries}\label{sec:preliminary}
\subsection{Primary decomposition theorem and the Jordan canonical form}\label{sec:aux}
In this subsection we recall known results that will be used in \cref{sec:Main}.
Recall that the minimal polynomial $q(x)$ of a matrix $\mathbf M\in\mathbb F^{I\times I}$ is the polynomial  of  least degree over  $\mathbb F$ whose leading coefficient is $1$
and  such that  $q(\mathbf M)=\mathbf O$. 
It is well known that  the minimal polynomial does not depend of $\mathbb F$, is unique,  and that the set of its  zeros  coincides with the set of the eigenvalues of the matrix (in the case $\mathbb F=\mathbb R$ both sets can be empty, namely, when the minimal polynomial does not have real roots).
Recall also that a non-constant polynomial is irreducible over $\mathbb F$ if its coefficients belong to $\mathbb F$ and it cannot be factorized into the product of two non-constant polynomials with coefficients in $\mathbb F$. 
For instance, the minimal polynomials of the matrices
$$
\begin{bmatrix}
0&0\\
1&1
\end{bmatrix},
\quad
\begin{bmatrix}
0&1\\
1&0
\end{bmatrix},
\quad
\begin{bmatrix}
0&1\\
-1&0
\end{bmatrix},
\quad
\begin{bmatrix}
0&1\\
0&0
\end{bmatrix},\text{ and } \mathbf I_I
$$
are $x^2-x$, $x^2-1$, $x^2+1$, $x^2$, and $x-1$, respectively.  The matrix $\mathbf I_I$ has a single eigenvalue $ 1$ of  multiplicity $I$ which corresponds to a  single root of $x-1$ of  multiplicity $1$. The polynomial $x^2+1$ is irreducible over $\mathbb R$ and is reducible over $\mathbb C$, $x^2+1=(x+i)(x-i)$, which agrees with the fact that the matrix $\begin{bmatrix}
0&1\\
-1&0
\end{bmatrix}$ does not have eigenvalues over $\mathbb R$ but has two eigenvalues $-i$ and $i$ over $\mathbb C$. 
It is well known that any polynomial with leading coefficient $1$ can be factorized as 
$$
q(x)=p_1(x)^{\mu_1}\cdots p_R(x)^{\mu_R}
$$
where  $p_r$ are distinct irreducible polynomials and $\mu_r\geq 1$.  Since in this paper $\fF$ is either $\mathbb C$ or $\mathbb R$, we have that 
\begin{align*}
&p_1,\dots,p_R\in\{x-\lambda:\ \lambda\in \mathbb C\}, \text{ if } \mathbb F=\mathbb C, \\
&p_1,\dots,p_R\in\{x-\lambda:\ \lambda\in \mathbb R\}\cup\{x^2+2ax+a^2+b^2:\ a,b\in \mathbb R \text{ and } b>0\},\text{ if }\ \mathbb F=\mathbb R.
\end{align*}
 The following theorem implies that the minimal polynomial of a matrix can be used to construct a basis in  which that matrix has block-diagonal form. 
\begin{theorem}[Primary decomposition theorem {\cite[pp.196--197]{Curtis_1984}}]\label{thm:PDT}
Let $\mathbf M\in\mathbb F^{I\times I}$ and let 
$$
q(x)=p_1(x)^{\mu_1}\cdots p_R(x)^{\mu_R}
$$
be the minimal polynomial of $\mathbf M$, factorized into powers of distinct  polynomials $p_r(x)$ that are irreducible (over $\mathbb F$). Then the subspaces
$$
E_r:=\nullsp{p_r(\mathbf M)^{\mu_r}},\quad  1\leq r\leq R
$$
 are invariant for $\mathbf M$,  i.e., $\mathbf M E_r\subseteq E_r$ and we have
\begin{equation}
\mathbb F^I =E_1\oplus\cdots\oplus E_R,\label{eq:dirsumdecomp}
\end{equation}
where ``$\oplus$'' denotes the direct sum of subspaces.
\end{theorem}
Decomposition  \cref{eq:dirsumdecomp} in \cref{thm:PDT} implies that  the matrix $\mathbf M$ is similar to a block-diagonal matrix.
Indeed, let $L_r=\dim E_r$ and let the columns of $\mathbf S_r\in \mathbb F^{I\times L_r}$ form a basis of $E_r$, $r=1,\dots, R$.  Then by \cref{eq:dirsumdecomp}, the columns of $\mathbf S:=[\mathbf S_1\ \dots\ \mathbf S_R]$ form a basis of
the entire space $\mathbb F^I$, implying that $\mathbf S$ is nonsingular.  Since $\mathbf M E_r \subseteq E_r$ it follows that there exists a unique matrix 
$\mathbf T_r\in\mathbb F^{L_r\times L_r}$ such that $\mathbf M \mathbf S_r = \mathbf S_r \mathbf T_r$, $r=1,\dots, R$. Hence  $\mathbf M[\mathbf S_1\ \dots\ \mathbf S_R] = [\mathbf S_1\mathbf T_1\ \dots\ \mathbf S_R\mathbf T_R]$ or
$$
\mathbf M=\mathbf S\Bdiag(\mathbf T_{1},\dots,\mathbf T_{R})\mathbf S^{-1}, \quad \mathbf S=[\mathbf S_1\ \dots\ \mathbf S_R], \quad
\mathbf S_r\in\mathbb F^{I\times L_r}.
$$
It is well-known that each of the matrices $\mathbf T_{r}$ can   further  be reduced  to  Jordan canonical form by a similarity transform.
Namely, if $p_r(x)^{\mu_r}=(x-\lambda)^{\mu_r}$ with $\lambda\in\mathbb F$, then
$\mathbf T_{r}$ is similar to $J(\lambda,n_{r1})\oplus\dots\oplus J(\lambda,n_{rk_r})$, where
$J(\lambda,n)$ denotes the $n\times n$ Jordan block with $\lambda$ on the main diagonal:
$$
\begin{bmatrix}
\lambda&       1&       0&\dots&0\\
0      & \lambda&       1&\dots&0\\
\vdots & \vdots &\vdots  &\vdots&\vdots\\
0      & 0      &       0&\dots&1\\
0      & 0      &       0&\dots&\lambda
\end{bmatrix}.
$$
 If $\mathbb F=\mathbb R$ and 
$p_r(x)^{\mu_r}=(x^2+2ax+a^2+b^2)^{\mu_r}$ with $a,b\in\mathbb R$ and $b>0$, then 
$\mathbf T_{r}$ is similar to $C(a,b,n_{r1})\oplus\dots\oplus C(a,b,n_{rk_r})$, where
$C(a,b,n)$ denotes the $2n\times 2n$  block matrix of the form
$$
\begin{bmatrix}
C(a,b)&       \mathbf I_2&       0&\dots&0\\
0      & C(a,b)&       \mathbf I_2&\dots&0\\
\vdots & \vdots &\vdots  &\vdots&\vdots\\
0      & 0      &       0&\dots&\mathbf I_2\\
0      & 0      &       0&\dots&C(a,b)
\end{bmatrix},\qquad 
C(a,b)=
\begin{bmatrix}
a&b\\
-b&a
\end{bmatrix}.
$$
It is  known that the values $n_{r_1},\dots,n_{rk_r}$ are uniquely determined by $\mathbf T_{r}$ up to permutation, in particular,
	$\max(n_{r_1},\dots,n_{rk_r})=\mu_r$. Thus, the Jordan canonical form is unique up to permutation of its  blocks.
 For more details on the Jordan canonical form we refer to \cite[Chapter 3]{HornJohnson}.
\subsection{An auxiliary result about simultaneous compression of a pair of tensors}\label{secmovedlemma}
Let ${\mathcal A},{\mathcal B}\in\fF^{{I}_1\times\dots\times {I}_N}$.
It is clear that the conditions
\begin{equation}
\operatorname{col}(\unfgen{{\mathbf B}}{n}) \subseteq \operatorname{col}(\unfgen{{\mathbf A}}{n}), \quad n\in\{1,\dots, \upN\}.
\label{eq:colAn=colBnnotilde}
\end{equation}
can be rewritten as
\begin{equation}
\unfgen{{\mathbf B}}{n}= \unfgen{{\mathbf A}}{n}\mathbf M_n,\qquad n\in\{1,\dots,\upN\},
\label{eq:matricesM_nmotilde}
\end{equation}
in which  $\mathbf M_n\in\mathbb F^{I_n\times I_n}$ is not necessarily unique. The goal of the following lemma is to show that  \cref{eq:matricesM_nmotilde} can further be reduced to the case where the matrices $\unfgen{{\mathbf A}}{n}$ do have full column rank, so $\mathbf M_n$ can be uniquely recovered as $\mathbf M_n={{\mathbf A}^\dagger_{(n^c;n)}}\unf{B}{n}$. In \cref{sec:Mainconnections} we will use $\mathbf M_1,\dots,\mathbf M_{\upN}$ to	 establish
connections between the terms in the decompositions of $\mathcal A$ and $\mathcal B$.
\begin{lemma}\label{lemma:moved}    
	Let $\tilde{\mathcal A},\tilde{\mathcal B}\in\fF^{\tilde{I}_1\times\dots\times \tilde{I}_N}$, $N\geq\upN\geq 2$ and let $\tilde{\mathcal A}$ be ML rank-$(I_1,\dots,I_{\upN},\cdot,\dots,\cdot)$. Assume that 
	\begin{equation}
	\operatorname{col}(\unfgen{\tilde{\mathbf B}}{n}) \subseteq \operatorname{col}(\unfgen{\tilde{\mathbf A}}{n}), \quad n\in\{1,\dots, \upN\}.
	\label{eq:colAn=colBntilde}
	\end{equation}
		Let also the rows of 
	$\mathbf U_n\in\fF^{I_n\times \tilde{I}_n}$ form an orthonormal basis of the row space of $\unfgen{\tilde{\mathbf A}}{n}$, $n\in\{1,\dots,\upN\}$\footnote{{For instance, one can take $\mathbf U_n$ equal to the transpose of the  ``$U$'' factor in the compact SVD of 	{$\tilde{\mathbf A}^T_{(n^c;n)}$}.} In this case, \cref{eq:tildeA=AtildeB=B} implements a standard compression by multilinear singular value decomposition \cite{2000LDLetall, Tucker1964, Tucker1966}, in which the compression matrices are obtained from $\mathcal A$.} and
	\begin{equation}
	\label{eq:A=tildeAB=tildeB}
	\mathcal A:=\tilde{\mathcal A}\bullet_1\mathbf U_1^*\cdots \bullet_{\upN}\mathbf U_{\upN}^*,\qquad 
	\mathcal B:=\tilde{\mathcal B}\bullet_1\mathbf U_1^*\cdots \bullet_{\upN}\mathbf U_{\upN}^*.
	\end{equation}
	Then the following statements hold.
	\begin{enumerate}
		\item \label{item:auxlemma0}
			For all $k\in\{1,\dots,N\}$, the row space of
			$\unfgen{\tilde{\mathbf A}}{k}$ contains the row space of $\unfgen{\tilde{\mathbf B}}{k}$. 
		\item \label{item:auxlemma2}
		$\tilde{\mathcal A}$ and $\tilde{\mathcal B}$ can be recovered from $\mathcal A$ and $\mathcal B$, respectively, as
		\begin{equation}
		\label{eq:tildeA=AtildeB=B}
		\tilde{\mathcal A}=\mathcal A\bullet_1\mathbf U_1^T\cdots \bullet_{\upN}\mathbf U_{\upN}^T,\qquad 
		\tilde{\mathcal B}=\mathcal B\bullet_1\mathbf U_1^T\cdots \bullet_{\upN}\mathbf U_{\upN}^T.
		\end{equation}
		\item \label{item:auxlemma1}
		$\mathcal A,\mathcal B\in\fF^{I_1\times\dots\times I_{\upN}\times \tilde{I}_{\upN+1}\times\dots\times \tilde{I}_{N}}$,
		$\mathcal A$ is ML rank-$(I_1,\dots,I_{\upN},\cdot,\dots,\cdot)$, and the
		ML rank of $\mathcal B$ equals the ML rank of $\tilde{\mathcal B}$.
	\end{enumerate}
\end{lemma}  
\begin{proof}
	\ref{item:auxlemma0}. Recall that  \cref{eq:An_one} is equivalent to any identity in \cref{eq:Ak_remaining}. Hence if
		\cref{eq:An_one} holds for $n=1$ and $n=2$, then, by \cref{eq:Ak_remaining},  the row space of
		$\unfgen{\mathbf D}{k}$ contains the row space of $\unfgen{\mathbf A}{k}$ for $k\in\{2,\dots,N\}$ and for $k\in\{1,3,\dots,N\}$, respectively, i.e., for
		all $k$. To complete the proof one should replace $\mathcal D$ and $\mathcal A$ in  \cref{eq:An_one,eq:Ak_remaining} by 
		$\tilde{\mathcal A}$ and $\tilde{\mathcal B}$, respectively.
	
	\ref{item:auxlemma2}.
	Since the  rows of $\mathbf U_n$ form an orthonormal basis of the row space of $\unfgen{\tilde{\mathbf A}}{n}$, it follows that
	$\unfgen{\tilde{\mathbf A}}{n}\mathbf U_n^H\mathbf U_n = \unfgen{\tilde{\mathbf A}}{n}$ or $\tilde{\mathcal A}\bullet_n (\mathbf U_n^T\mathbf U_n^*)=\tilde{\mathcal A}$,
	$n\in\{1,\dots,\upN\}$.
	Hence 
	\begin{multline*}
	\mathcal A\bullet_1\mathbf U_1^T\cdots \bullet_{\upN}\mathbf U_{\upN}^T = 
	(\tilde{\mathcal A}\bullet_1\mathbf U_1^*\cdots \bullet_{\upN}\mathbf U_{\upN}^*)\bullet_1\mathbf U_1^T\cdots \bullet_{\upN}\mathbf U_{\upN}^T=\\
	\tilde{\mathcal A}\bullet_1(\mathbf U_1^T\mathbf U_1^*)\cdots \bullet_{\upN}(\mathbf U_{\upN}^T\mathbf U_{\upN}^*)=\tilde{\mathcal A}.
	\end{multline*}
	By statement \ref{item:auxlemma0}, the identity for $\tilde{\mathcal B}$ can be proved in a similar way.
	
	\ref{item:auxlemma1}.  From \cref{eq:An_one,eq:A=tildeAB=tildeB,eq:tildeA=AtildeB=B} it follows that
	$$
	r_{\unf{A}{n}}\leq r_{\unfgen{\tilde{\mathbf A}}{n}}\leq r_{\unf{A}{n}},\quad	
	r_{\unf{B}{n}}\leq r_{\unfgen{\tilde{\mathbf B}}{n}}\leq r_{\unf{B}{n}},\quad n=1,\dots,\hat N
	$$
	implying that $r_{\unf{A}{n}}= r_{\unfgen{\tilde{\mathbf A}}{n}}=I_n$ and 	$r_{\unf{B}{n}}=r_{\unfgen{\tilde{\mathbf B}}{n}}$ for $n=1,\dots,\hat N$.
	\end{proof}

\section{Main results}\label{sec:Main}
\subsection{Connections between tensors $\mathcal A$ and $\mathcal B$  that satisfy the first $\upN$  conditions in \cref{eq:AincludesBnew}}\label{sec:Mainconnections}
To simplify the presentation throughout this subsection  we assume that 
the first $\upN$ matrix  representations of $\mathcal A$ have full column rank. The general case follows from  \cref{lemma:moved} above.
Also, to keep the presentation and derivation of results easy to follow,  we first consider the particular case where  
$\mathcal A$ and $\mathcal B$ are third-order  tensors (i.e., $N=3$) that satisfy only
the first two conditions (i.e., $\upN=2$) in 
\begin{equation}
\operatorname{col}(\unfone{B})\subseteq \operatorname{col}(\unfone{A}),\ 
\operatorname{col}(\unftwo{B})\subseteq \operatorname{col}(\unftwo{A}),\ 
\operatorname{col}(\unfthree{B})\subseteq \operatorname{col}(\unfthree{A}). \label{eq:three inclusions}
\end{equation}
The case where all three conditions in \cref{eq:three inclusions} hold (i.e., $N=\upN=3$) and  the general case
$N\geq 3$, $N\geq \upN\geq 2$ 
will be covered  by \cref{thm:newverygeneral theorem} below. 
\begin{theorem}\label{thm3rdorder}
	Let tensors $\mathcal A, \mathcal B\in\mathbb F^{I_1\times I_2\times I_3}$. Assume that   
	\begin{equation}
	\unfone{A} \text{ and }\ \unftwo{A} \ \text{ have full column rank}\label{eq:X1X2fcr}
	\end{equation}
	and that there exist
 matrices $\mathbf M_1\in\mathbb F^{I_1\times I_1}$ and $\mathbf M_2\in\mathbb F^{I_2\times I_2}$ such that
	\begin{equation}
\unfone{B}= \unfone{A}\mathbf M_1\ \text{ and }\ 	\unftwo{B} = \unftwo{A}\mathbf M_2.\label{eq:A=BM}
	\end{equation}
	Then the following statements hold.
	\begin{enumerate}
		\item\label{item:1}
		The matrices $\mathbf M_1$ and $\mathbf M_2$ have the same minimal polynomial $q(x)$.
		\item\label{item:3}Consider the factorization $q(x)=p_1(x)^{\mu_1}\cdots p_R(x)^{\mu_R}$ with
distinct polynomials $p_r(x)$ that are irreducible (over $\mathbb F$) and set
$$
L_{1r}:=\dim(\nullsp{p_r(\mathbf M_1)^{\mu_r}}),\ L_{2r}:=\dim(\nullsp{p_r(\mathbf M_2)^{\mu_r}})\quad  1\leq r\leq R.
$$	
Let also
		\begin{align}
		\mathbf M_1&=\mathbf S_1\Bdiag(\mathbf T_{11},\dots,\mathbf T_{1R})\mathbf S_1^{-1}, \quad \mathbf S_1=[\mathbf S_{11}\ \dots\ \mathbf S_{1R}], \quad
		\mathbf S_{1r}\in\mathbb F^{I_1\times L_{1r}},\label{eq:PMD1}\\
		\mathbf M_2&=\mathbf S_2\Bdiag(\mathbf T_{21},\dots,\mathbf T_{2R})\mathbf S_2^{-1}, \quad  \mathbf S_2=[\mathbf S_{21}\ \dots\ \mathbf S_{2R}], \quad
				\mathbf S_{2r}\in\mathbb F^{I_2\times L_{2r}}\label{eq:PMD2}
		\end{align}
		be the primary decompositions of $\mathbf M_1$ and $\mathbf M_2$, respectively,
		such that the minimal polynomials of $\mathbf T_{1r}$ and $\mathbf T_{2r}$   are equal to $p_r(x)^{\mu_r}$ for each $r=1,\dots,R$.
		Then the matrices 
		\begin{equation}
		\mathbf D_i:=\mathbf S_1^T\mathbf A_i\mathbf S_2,\qquad \mathbf A_i := \mathbf A(:,:,i),\qquad  i=1,\dots,I_3
		\label{eq:DeqSAS}
		\end{equation} are block-diagonal,
		$
		\mathbf D_i=\Bdiag(\mathbf D_{i,11},\dots,\mathbf D_{i,RR})$, $\mathbf D_{i,rr}\in\mathbb F^{L_{1r}\times L_{2r}}
		$ and 
		\begin{equation}
		\mathbf T_{1r}^T\mathbf D_{i,rr}=\mathbf D_{i,rr}\mathbf T_{2r},\qquad i=1,\dots,I_3,\quad r=1,\dots,R.
		\label{eq:TDDTfirst}
		\end{equation}
		\item\label{tem:shouldbe3} Let $\mathcal D_r\in\mathbb F^{L_{1r}\times L_{2r}\times I_3}$ denote a tensor with frontal slices
		$\mathbf D_{1,rr},\dots,\mathbf D_{I_3,rr}\in\mathbb F^{L_{1r}\times L_{2r}}$ and let 
		\begin{align*}
		\mathbf S_1^{-T} =: \mathbf X^{(1)}= [\mathbf X^{(1)}_1\ \dots\ \mathbf X^{(1)}_R],\quad \mathbf X^{(1)}_r\in\mathbb F^{I_1\times L_{1r}},\\
   		\mathbf S_2^{-T} =: \mathbf X^{(2)}= [\mathbf X^{(2)}_1\ \dots\ \mathbf X^{(2)}_R],\quad \mathbf X^{(2)}_r\in\mathbb F^{I_2\times L_{2r}}.
		\end{align*}
		Then
		 the tensors $\mathcal A$ and $\mathcal B$ admit  decompositions into  ML rank-$(L_{1r}, L_{2r},\cdot)$ terms which are connected as follows: 
		 \begin{align}
		 		\mathcal A &= \sum\limits_{r=1}^R \mathcal D_r \bullet_1 \mathbf X^{(1)}_r \bullet_2\mathbf X^{(2)}_r=:
		 		\sum\limits_{r=1}^R\mathcal A_{r},\label{eq:decA1r}\\
		 		\mathcal B &= \sum\limits_{r=1}^R (\mathcal D_r\bullet_1 \mathbf T_{1r}^T)\bullet_1 \mathbf X^{(1)}_r \bullet_2\mathbf X^{(2)}_r=:
		 		\sum\limits_{r=1}^R\mathcal B_{r},
		 		\label{eq:decofB11}
		 \end{align}		
		 and
		 \begin{equation}
		 \mathcal D_r\bullet_1 \mathbf T_{1r}^T=\mathcal D_r\bullet_2\mathbf T_{2r}^T\qquad r=1,\dots,R.
		 \label{eq:TDDT}
		 \end{equation}
        \item  \label{item:4}
		If $I_1=I_2$ and if there exists a  linear combination of $\mathbf A_1,\dots,\mathbf A_{I_3}$ that is nonsingular, then 
		$\mathbf M_1$ is similar to $\mathbf M_2$.
		\item \label{item:4part2}  If $\mathbf M_1$ is similar to $\mathbf M_2$, then 
		$L_{1r}=L_{2r}$ for all $r$ and the matrices $\mathbf S_1$ and $\mathbf S_2$ in \cref{eq:PMD1} and \cref{eq:PMD2} can be chosen such that
		$\mathbf T_{1r}=\mathbf T_{2r}$ for all $r$.
		\item \label{item:5}
		 If, for some $r$, the matrix $\mathbf T_{1r}$  (or $\mathbf T_{2r}$) is a scalar multiple of the identity matrix, i.e., if  $\mathbf T_{1r}=\lambda_r\mathbf I_{L_{1r}}$  (or $\mathbf T_{2r}=\lambda_r\mathbf I_{L_{2r}}$), then $\mathcal A_{r}=\lambda_r \mathcal B_r$. 
		\item \label{item:6}
		 If    $\mathbf T_{1r}=\lambda_r\mathbf I_{L_{1r}}$  (or $\mathbf T_{2r}=\lambda_r\mathbf I_{L_{2r}}$)  for all $r$, then 
		 $\mathcal A$ and $\mathcal B$ consist of the same ML rank-$(L_{1r},L_{2 r},\cdot)$ terms, possibly differently scaled.
	\end{enumerate}
\end{theorem}
\begin{proof}
	\ref{item:1}.	To prove that the minimal polynomials of $\mathbf M_1$ and $\mathbf M_2$ coincide, it is sufficient to show that
	a polynomial $q(x)$ annihilates $\mathbf M_1$ if and only if $q(x)$ annihilates $\mathbf M_2$.
  	By	\cref{eq:A=BM}, $\mathcal B=\mathcal A\bullet_1 \mathbf M_1^T=\mathcal A\bullet_2 \mathbf M_2^T$.
	Since, by \cref{eq:mode_n-unfolding},
	\begin{equation}
	\unfone{A}=[\mathbf A_1\ \dots\ \mathbf A_{I_3}]^T\text{ and } \unftwo{A}=[\mathbf A_1^T\ \dots\ \mathbf A_{I_3}^T]^T,
	\label{eq:B1B2unf}
	\end{equation}
	it follows that
	\begin{equation}
	(\mathbf B_i=)\mathbf M_1^T\mathbf A_i = \mathbf A_i\mathbf M_2,\qquad i\in\{1,\dots,I_3\}.\label{eq:mainidentity}
	\end{equation}
	Hence for any $k\geq 1$, 
	\begin{align*}
	(\mathbf M_1^T)^k\mathbf A_i =& (\mathbf M_1^T)^{k-1}\mathbf M_1^T\mathbf A_i    = (\mathbf M_1^T)^{k-1}  \mathbf A_i\mathbf M_2=\\
	&(\mathbf M_1^T)^{k-2}\mathbf M_1^T  \mathbf A_i\mathbf M_2=
	(\mathbf M_1^T)^{k-2}  \mathbf A_i\mathbf M_2^2=\dots= \mathbf A_i\mathbf M_2^k,
	\end{align*}
	implying that for any polynomial $q$,
	\begin{equation}
	q(\mathbf M_1)^T\mathbf A_i = \mathbf A_i q(\mathbf M_2),\qquad i\in\{1,\dots,I_3\}.
	\label{eq:minpol}
	\end{equation}
	It follows from \cref{eq:B1B2unf} that  \cref{eq:minpol} is equivalent to 
	\begin{equation}
       \unftwo{A}q(\mathbf M_2) = \Bdiag(q(\mathbf M_1)^T,\dots,q(\mathbf M_1)^T)\unftwo{A}
       \label{eq:eqb2q}
	\end{equation}
	and to
	\begin{equation}
       \unfone{A}q(\mathbf M_1) = \Bdiag(q(\mathbf M_2)^T,\dots,q(\mathbf M_2)^T)\unfone{A}.
       \label{eq:eqb1q}
	\end{equation}
	Assume that $q$ annihilates $\mathbf M_1$. Then, by \cref{eq:eqb2q}, $\unftwo{A}q(\mathbf M_2)=\mathbf O$. Since
	$\unftwo{A}$ has full column rank, it follows that $q$ annihilates $\mathbf M_2$. On the other hand, if $q$ annihilates $\mathbf M_2$, then
	by \cref{eq:eqb1q}, $\unfone{A}q(\mathbf M_1)=\mathbf O$.  Since $\unfone{A}$ has full column rank, it follows that $q$ annihilates $\mathbf M_1$. Thus, the matrices $\mathbf M_1$ and $\mathbf M_2$ have the same minimal polynomial.
	
	\ref{item:3}.  By \cref{eq:PMD1,eq:PMD2,eq:mainidentity},
	\begin{multline*}
	(\mathbf S_1\Bdiag(\mathbf T_{11},\dots,\mathbf T_{1R})\mathbf S_1^{-1})^T\cdot \mathbf A_i
	= \\
	\mathbf A_i\cdot  \mathbf S_2\Bdiag(\mathbf T_{21},\dots,\mathbf T_{2R})\mathbf S_2^{-1}
	,\qquad i\in\{1,\dots,I_3\}.
	\end{multline*}
	Hence
	\begin{multline}
	\Bdiag(\mathbf T_{11}^T,\dots,\mathbf T_{1R}^T)
	\mathbf S_1^T\mathbf A_i\mathbf S_2 =\\
	\mathbf S_1^T\mathbf A_i\mathbf S_2
	\Bdiag(\mathbf T_{21},\dots,\mathbf T_{2R}),\qquad i\in\{1,\dots,I_3\}.
	\label{eq:mltline}
	\end{multline}
	Let 
	$$
	\mathbf S_1^T\mathbf A_i\mathbf S_2=:\mathbf D_i=\left(\mathbf D_{i,r_1r_2}\right)_{r_1,r_2=1}^R
	$$
	 denote a block matrix with
	$\mathbf D_{i,r_1r_2}\in\mathbb F^{L_{1r_1}\times L_{2r_2}}$.
	It is clear that  \cref{eq:mltline} can be rewritten as
	\begin{equation}
	\mathbf T_{1r_1}^T \mathbf D_{i,r_1r_2} = \mathbf D_{i,r_1r_2} \mathbf T_{2r_2},\qquad r_1,r_2=1,\dots,R,\qquad i\in\{1,\dots,I_3\},
	\label{eq:sylv}
	\end{equation}
	implying that \cref{eq:TDDTfirst} holds.
	
	Now we show that $\mathbf D_i$ is a block diagonal matrix, i.e., that $\mathbf D_{i,r_1r_2}=\mathbf O$ for $r_1\ne r_2$. 
	Let $p_r(x)^{\mu_r}$ denote the minimal polynomial of $\mathbf T_{1r}$ (or $\mathbf T_{2r}$). Then, by \cref{eq:sylv},   
	\begin{equation}
	\mathbf O=\left(p_{r_1}(\mathbf T_{1r_1})^{\mu_{r_1}}\right)^T \mathbf D_{i,r_1r_2} = \mathbf D_{i,r_1r_2} p_{r_1}(\mathbf T_{2r_2})^{\mu_{r_1}},
	\label{eq:sylvwithpoly}
	\end{equation}	
	for all  $r_1,r_2=1,\dots,R$ and $i\in\{1,\dots,I_3\}$. Let $r_1\ne r_2$. To prove that $\mathbf D_{i,r_1r_2}=\mathbf O$, it is sufficient to show that the matrix  $p_{r_1}(\mathbf T_{2r_2})^{\mu_{r_1}}$  is nonsingular. Since the polynomials $p_{r_1}(x)^{\mu_{r_1}}$ 
	and $p_{r_2}(x)^{\mu_{r_2}}$ are relatively prime, it follows  from the Euclidean algorithm that there exist polynomials $f(x)$ and $g(x)$ such that
	$1=p_{r_1}(x)^{\mu_{r_1}}f(x)+p_{r_2}(x)^{\mu_{r_2}}g(x)$ for all $x\in\mathbb F$. Hence
	$$
	\mathbf I = p_{r_1}(\mathbf T_{2r_2})^{\mu_{r_1}}f(\mathbf T_{2r_2})+p_{r_2}(\mathbf T_{2r_2})^{\mu_{r_2}}g(\mathbf T_{2r_2})=
	p_{r_1}(\mathbf T_{2r_2})^{\mu_{r_1}}f(\mathbf T_{2r_2}).
	$$
	Thus, $p_{r_1}(\mathbf T_{2r_2})^{\mu_{r_1}}$  is nonsingular.

   \ref{tem:shouldbe3}.	By \cref{eq:DeqSAS},
       \begin{equation}
   		\mathbf A_i=\mathbf S_1^{-T}\mathbf D_i\mathbf S_2^{-1}=\mathbf X^{(1)}\mathbf D_i\mathbf X^{(2)T},\qquad  i=1,\dots,I_3
   		\label{eq:AeqXDX}
   		\end{equation}
   		which is equivalent to \cref{eq:decA1r}.
Since, by \cref{eq:A=BM}, $\mathbf B_i=\mathbf M_1^{T}\mathbf A_i$, it follows from \cref{eq:PMD1,eq:AeqXDX} that
        \begin{multline}
        \mathbf B_i=\mathbf M_1^{T}\mathbf A_i =(\mathbf S_1\Bdiag(\mathbf T_{11},\dots,\mathbf T_{1R})\mathbf S_1^{-1})^{T} \mathbf S_1^{-T}\mathbf D_i\mathbf S_2^{-1}\\
   		\mathbf S_1^{-T}\Bdiag(\mathbf T_{11}^{T},\dots,\mathbf T_{1R}^{T})\mathbf S_1^{T}
   		\mathbf S_1^{-T}\mathbf D_i\mathbf S_2^{-1}=\\
   		\mathbf X^{(1)}\Bdiag(\mathbf T_{11}^{T}\mathbf D_{i,11},\dots,\mathbf T_{1R}^{T}\mathbf D_{i,RR})\mathbf X^{(2)T},\qquad  i=1,\dots,I_3,
   		\label{eq:BieqXDX}
   		\end{multline}
 which is equivalent to \cref{eq:decofB11}.  Finally, identity \cref{eq:TDDT}	is equivalent to \cref{eq:TDDTfirst}.	
	
	\ref{item:4}.  Let the linear combination $t_1\mathbf A_1+\dots+t_{I_3}\mathbf A_{I_3}$ be nonsingular. Then, by \cref{eq:mainidentity},
	$$
	\mathbf M_2 = (t_1\mathbf A_1+\dots+t_{I_3}\mathbf A_{I_3})^{-1}\mathbf M_1^T (t_1\mathbf A_1+\dots+t_{I_3}\mathbf A_{I_3}),
	$$
	i.e., $\mathbf M_2$ is similar to $\mathbf M_1^T$.  Since any matrix is similar to its transpose \cite[Section 3.2.3]{HornJohnson}, it follows that
	$\mathbf M_2$ is similar to $\mathbf M_1$.
	
	\ref{item:4part2}.  We choose $\mathbf S_1$ such that  the matrices  $\mathbf T_{11}, \dots, \mathbf T_{1R}$ in \cref{eq:PMD1} are in the Jordan canonical form.  Since similar matrices have the same Jordan canonical form,
	the matrix $\mathbf M_2$ is similar to $\Bdiag(\mathbf T_{11},\dots,\mathbf T_{1R})$, i.e., there exists $\mathbf S_2$ such that \cref{eq:PMD2} holds for $\mathbf T_{11}=\mathbf T_{21},\dots,  \mathbf T_{1R}=\mathbf T_{2R}$.
	
	\ref{item:5}. and \ref{item:6}. follow from \cref{eq:decofB11}. 	
\end{proof}
\begin{expl}
   This example illustrates that although the matrices $\mathbf M_1$ and $\mathbf M_2$  in   \cref{thm3rdorder} have the same minimal polynomial they are not necessarily similar.
	Let the frontal slices  of $\mathcal A\in\mathbb C^{3\times 3\times 4}$ have the following nonzero pattern:
	$$
	\begin{bmatrix}
	0&*&*\\ *&0&0\\ *&0&0
	\end{bmatrix}
	$$
	It is clear that any linear combination of the frontal slices of $\mathcal A$ is singular so the assumption in statement \ref{item:4} of \cref{thm3rdorder} does not hold.
	We choose the values ``$*$''  (e.g., generic values) such that   $\unfone{A}$ and $\unftwo{A}$ have full column rank.  It is clear that $\mathcal A$ is the sum of a ML rank-$(1,2,\cdot)$ and a ML rank-$(2,1,\cdot)$ term. More precisely, $\mathcal A$ 	is the sum of a ML rank-$(1,2,2)$ and a ML rank-$(2,1,2)$ term. 
	Let $\mathbf M_1:=\operatorname{diag}(\lambda_2,\lambda_1,\lambda_1)$ and $\mathcal B=\mathcal A\bullet_1 \mathbf M_1^T$.
	One can easily verify that $\mathcal B=\mathcal A\bullet_2 \mathbf M_2^T$, where
	$\mathbf M_2=\operatorname{diag}(\lambda_1,\lambda_2,\lambda_2)$.
	Thus, if $\lambda_1\ne \lambda_2$, then $\mathbf M_1$ and $\mathbf M_2$ have the same minimal polynomial but are not similar. 
\end{expl}

Now we consider the general case, that is, we assume that $\mathcal A$ and $\mathcal B$ are tensors of order $N\geq 3$ and satisfy \cref{eq:colAn=colBnnotilde} for $N\geq\upN\geq 2$. First we extend the notion of  block diagonal matrices to tensors.
Let the numbers $L_{n1},\dots,L_{nR}$ sum up to $I_n$ for each $n=1,\dots,\upN$. Consider the partition of $\{1,\dots,I_n\}$ into consecutive blocks $V_{n1},\dots,V_{nR} $ of length $L_{n1},\dots,L_{nR}$, respectively, so $V_{n1}=\{1,\dots,L_{n1}\},\dots, V_{nR}=\{I_n-L_{nR}+1,\dots,L_{nR}\}$. 
If the  condition
\begin{equation}
(\mathcal D)_{i_1,\dots,i_N}=0\text{ for } (i_1,\dots,i_{\upN})\not\in\bigcup\limits_{r=1}^R (V_{1r}\times\dots\times V_{\upN r})
\label{eq:blockdiagonaltens}
\end{equation}
holds, then we say that  $\mathcal D$ is a block diagonal tensor and write $\mathcal D=\Bdiag(\mathcal D_{1},\dots,\mathcal D_{R})$,
where $\mathcal D_{r}
:= \mathcal D(V_{1r},\dots,V_{\upN r},:,\dots,:)\in\mathbb F^{L_{1r}\times\dots\times L_{\upN r}\times I_{\upN+1}\times\dots\times I_N }$ denote the diagonal blocks. For instance, 
statement \ref{item:3} of \cref{thm3rdorder} means that
if $\mathcal D$ is the $I_1\times I_2\times I_3$ tensor formed by  the $I_1\times I_2$  matrices  $\mathbf D_i$  in \cref{eq:DeqSAS}, i.e., if 
$\mathcal D:=\mathcal A  \bullet_1 \mathbf S_1 \bullet_{2} \mathbf S_{2}$, then $\mathcal D=\Bdiag(\mathcal D_{1},\dots,\mathcal D_{R})$, where
the diagonal blocks $\mathcal D_{r}\in\mathbb F^{L_{1r}\times L_{2 r}\times I_{3}}$ are defined in statement \ref{tem:shouldbe3} of  \cref{thm3rdorder}.

The following result generalizes  \cref{thm3rdorder} for $N\geq 3$ and  $N\geq \upN\geq 2$.
\begin{theorem}\label{thm:newverygeneral theorem}
	Let tensors $\mathcal A, \mathcal B\in\mathbb F^{I_1\times \dots\times I_N}$ and let $N\geq 3$, $N\geq\upN\geq 2$. Assume that   for each 
	$n\in\{1\dots,\upN\}$,
	\begin{equation}
	\unf{A}{n}  \text{ has full column rank} \label{eq:X1X2fcrHO}
	\end{equation}
	and that there exists
	matrix $\mathbf M_n\in\mathbb F^{I_n\times I_n} $ such that
	\begin{equation}
	\unf{B}{n}= \unf{A}{n}\mathbf M_n.\label{eq:A=BMHO}
	\end{equation}
	Then the following statements hold.
	\begin{enumerate}
		\item\label{item:1HO}
		The matrices $\mathbf M_1,\dots,\mathbf M_{\upN}$ have the same minimal polynomial $q(x)$.
		\item\label{item:3HO}Consider the factorization $q(x)=p_1(x)^{\mu_1}\cdots p_R(x)^{\mu_R}$ with
			distinct polynomials $p_r(x)$ that are irreducible (over $\mathbb F$) and set
		$$
		L_{nr}:=\dim(\nullsp{p_r(\mathbf M_{n})^{\mu_r}}),\qquad  1\leq r\leq R, \quad 1\leq n\leq\upN.
		$$
			Let also
		\begin{equation}
		\begin{split}
		\mathbf M_1&=\mathbf S_1\Bdiag(\mathbf T_{11},\dots,\mathbf T_{1R})\mathbf S_1^{-1}, \quad \mathbf S_1=[\mathbf S_{11}\ \dots\ \mathbf S_{1R}], \quad
		\mathbf S_{1r}\in\mathbb F^{I_1\times L_{1r}},
		\\
		&\phantom{=}\vdots
		\\
		\mathbf M_{\upN}&=\mathbf S_{\upN}\Bdiag(\mathbf T_{{\upN}1},\dots,\mathbf T_{{\upN}R})\mathbf S_{\upN}^{-1}, \quad  \mathbf S_{\upN}=[\mathbf S_{{\upN}1}\ \dots\ \mathbf S_{{\upN}R}], \quad
		\mathbf S_{{\upN}r}\in\mathbb F^{I_{\upN}\times L_{{\upN}r}}
		\end{split}
		\label{eq:PMD1HO}
		\end{equation}
		be the primary decompositions of $\mathbf M_1,\dots,\mathbf M_{\upN}$, respectively,
		such that the minimal polynomials of $\mathbf T_{1r},\dots,\mathbf T_{{\upN}r}$ are equal to $p_r(x)^{\mu_r}$ for each $r=1,\dots,R$. Then the tensor 
		$$
		\mathcal D:=\mathcal A  \bullet_1 \mathbf S_1\dots \bullet_{\upN} \mathbf S_{\upN}
		$$
		is block-diagonal (see \cref{eq:blockdiagonaltens}),
		$$
		\mathcal D=\Bdiag(\mathcal D_{1},\dots,\mathcal D_{R}),\qquad  \mathcal D_{r}\in\mathbb F^{L_{1r}\times\dots\times L_{\upN r}\times I_{\upN+1}\times\dots\times I_N }
		$$ 
         and 
       \begin{equation}
       \mathcal D_r\bullet_1 \mathbf T_{1r}^T=\dots=\mathcal D_r\bullet_{\upN}\mathbf T_{{\upN}r}^T\qquad r=1,\dots,R.
       \label{eq:TDDTHO}
       \end{equation}
		\item\label{tem:shouldbe3HO} 
		 Let 
		\begin{align*}
		\mathbf S_n^{-T} =: \mathbf X^{(n)}= [\mathbf X^{(n)}_1\ \dots\ \mathbf X^{(n)}_R],\quad \mathbf X^{(n)}_r\in\mathbb F^{I_n\times L_{nr}}.
		\end{align*}
		Then
		the tensors $\mathcal A$ and $\mathcal B$ admit  decompositions into  ML rank-$(L_{1r},\dots,L_{\upN r},\cdot,\dots,\cdot)$ terms which are connected as follows: 
		\begin{align}
		\mathcal A &= \sum\limits_{r=1}^R \mathcal D_r \bullet_1 \mathbf X^{(1)}_r\dots \bullet_{\upN}\mathbf X^{({\upN})}_r=:
		\sum\limits_{r=1}^R\mathcal A_{r},\label{eq:decA1rHO}\\
		\mathcal B &= \sum\limits_{r=1}^R (\mathcal D_r\bullet_1 \mathbf T_{1r}^T)\bullet_1 \mathbf X^{(1)}_r\dots \bullet_{\upN}\mathbf X^{({\upN})}_r=:
		\sum\limits_{r=1}^R\mathcal B_{r},
		\label{eq:decofB11HO}
		\end{align}		
		in which the tensors $\mathcal D_r$ satisfy the identities in \cref{eq:TDDTHO}.
		\item  \label{item:4HO}
		Let $\mathbf A_{ij,k}$, $k=1,\dots,(I_1\cdots I_N)/(I_iI_j)$ denote the $I_i\times I_j$ slices of $\mathcal A$, that is, 
		$\mathbf A_{ij,k}\in\mathbb F^{I_i\times I_j}$ is obtained from $\mathcal A$ by fixing all indices but $i$ and $j$.
		If $I_i=I_j$ and if there exists a  linear combination of $\mathbf A_{ij,k}$ that is nonsingular, then 
		$\mathbf M_i$ is similar to $\mathbf M_j$.
		\item \label{item:4part2HO}  If $\mathbf M_{i}$ is similar to $\mathbf M_{j}$, then 
		$L_{ir}=L_{jr}$ for all $r$ and the matrices $\mathbf S_{i}$ and $\mathbf S_{j}$ in \cref{eq:PMD1HO}  can be chosen such that
		$\mathbf T_{ir}=\mathbf T_{jr}$ for all $r$.
		\item \label{item:5HO}
		If, for some $r$, there exists $n$ such that  the matrix $\mathbf T_{nr}$   is a scalar multiple of the identity matrix, i.e., if   $\mathbf T_{nr}=\lambda_r\mathbf I_{L_{nr}}$, then $\mathcal A_{r}=\lambda_r \mathcal B_r$. 
		\item \label{item:6HO}
		If for each $r$  there exists $n_r$ such that $\mathbf T_{n_rr}=\lambda_r\mathbf I_{L_{n_rr}}$, then 
		$\mathcal A$ and $\mathcal B$ consist of the same ML rank-$(L_{1r},\dots,L_{\upN r},\cdot,\dots,\cdot)$ terms, possibly differently scaled.
	\end{enumerate}
\end{theorem}
\begin{proof}
Let $1\leq i<j\leq\upN$. We reshape 
$\mathcal A$ and $\mathcal B$ into the $I_i\times I_j\times \frac{\prod I_n}{I_{i}I_{j}}$ tensors 
$\mathcal A^{ij}$ and $\mathcal B^{ij}$ such that
\begin{equation}
{\mathbf A}^{ij}_{(2,3;1)}=\unfgen{\mathbf A}{i},\quad {\mathbf B}^{ij}_{(2,3;1)}=\unfgen{\mathbf B}{i},\quad
{\mathbf A}^{ij}_{(1,3;2)}=\unfgen{\mathbf A}{j},\quad {\mathbf B}^{ij}_{(1,3;2)}=\unfgen{\mathbf B}{j}.
\label{eq:redtoij}
\end{equation}
Then, by \cref{eq:redtoij,eq:X1X2fcrHO}, the first two matrix representations of  $\mathcal A^{ij}$ have full column rank and,
by \cref{eq:redtoij,eq:A=BMHO},
$$
{\mathbf B}^{ij}_{(2,3;1)}= {\mathbf A}^{ij}_{(2,3;1)}\mathbf M_{i} \text{ and } {\mathbf B}^{ij}_{(1,3;2)} = {\mathbf A}^{ij}_{(1,3;2)}\mathbf M_{j}.
$$
Thus $\mathcal A^{ij}$ and $\mathcal B^{ij}$ satisfy the assumptions in \cref{thm3rdorder}.
We leave it to the reader to show that the statements in \cref{thm:newverygeneral theorem}
can be obtained from the corresponding statements of \cref{thm3rdorder} by applying it to all pairs
($\mathcal A^{ij}$,$\mathcal B^{ij}$), where $1\leq i<j\leq\upN$.
\end{proof}
\subsection{Redundancy of  conditions in \cref{eq:manyinclusions}}\label{subsection:redundancy}
In this subsection we prove  that if  
$
\operatorname{col}(\unfset{\mathbf B}{n}) \subseteq \operatorname{col}(\unfset{\mathbf A}{n})
$,
then for any subset $S\subsetneq \{1,\dots,N\}$ that contains $n$ we also have  that
$
\operatorname{col}(\unfset{\mathbf B}{S })\\ \subseteq \operatorname{col}(\unfset{\mathbf A}{S})
$ (\cref{Thm:HOcase}).
 Hence the $N$ conditions in \cref{eq:AincludesBnew} imply the $2^N-2$ conditions in \cref{eq:manyinclusions} (\cref{corr:redandancy}).

Let us first formally define   generalized matrix representations. Let  $\mathcal A\in\fF^{I_1\times \dots \times I_N}$, let $S$ be a proper subset of  $\{1,\dots,N\}$  and let
$S^c$ denote the complement of $S$.
A {\em \moden{S} slice} of $\mathcal A$ is a  subtensor obtained from $\mathcal A$ by fixing the indices in $S$. It is clear that $\mathcal A$ has  $\prod\limits_{n\in S} I_n$ \moden{S} slices.
{\em A \moden{S} matrix  representation} of   $\mathcal A$ is a matrix  $\unfset{\mathbf A}{S}\in\fF^{(\prod\limits_{n\not\in S} I_n) \times (\prod\limits_{n\in S} I_n)}$  
whose columns are  the vectorized \moden{S} slices of $\mathcal A$.
Formally, if we follow the conventions that
\begin{equation}
S=\{q_1,\dots,q_{N-k}\}\text{ with } q_1<\dots<q_{N-k}\ \text{ and } \ S^c=\{p_1,\dots,p_k\}
\text{ with } p_1<\dots<p_k, 
\label{eq:SandSc}
\end{equation}
then
\begin{equation}
 \text{ the } (\colvec{p}{k},\colvec{q}{N-k})\text{th entry of  the matrix }\unfset{\mathbf A}{S}\text{  is equal to } a_{i_1\dots i_N},  
 \label{eq:defmatrreprS}
\end{equation}
where
$$
\colvec{p}{k} := 
1+\sum_{u=1}^k (i_{p_u}-1)\prod_{s=1}^{u-1}I_{p_s}
$$
denotes the linear index corresponding to the element in the $(i_{p_1},\dots, i_{p_k})$ position of an $I_{p_1}\times\dots\times I_{p_k}$ tensor.
If $S=\{n\}$, then $\unfset{\mathbf A}{S}$ coincides with the \moden{n} matrix representation $\unf{A}{n}$ introduced  earlier in \cref {eq:mode_n-unfolding}. 
In the following lemma  we prove that, if for two $I_1\times \dots \times I_N$  tensors $\mathcal A$ and $\mathcal B$ the identity  $ \unf{B}{n} = \unf{A}{n}\mathbf M_n$ holds for some $n$, then
for any subset $S$ that contains $n$ there exists a matrix $\mathbf M_S$ such that  $\unfset{\mathbf B}{S} = \unfset{\mathbf A}{S}\mathbf M_S$. In fact  equation \cref{eq:somepermutation} below implies  that the matrix
$\mathbf M_S$ coincides  up to column and row permutation with the direct sum of  $\mathbf M_n$  multiple times with itself. 
\begin{lemma}\label{Thm:HOcase}
	Let $N\geq 4$,  and let $\mathcal A, \mathcal B\in\mathbb F^{I_1\times \dots\times I_N}$ be  such that 
	$ \unf{B}{n} = \unf{A}{n}\mathbf M_n$ for some $\mathbf  M_n\in\mathbb F^{I_n\times I_n}$. Let  $S$ and $S^c$  be as in \cref{eq:SandSc} and let  $n\in S$, that is,  $q_l=n$ for some $l\in\{1,\dots,N-k\}$. Then $\unfset{\mathbf B}{S} = \unfset{\mathbf A}{S}{\mathbf M}_S$, where
	\begin{equation}
	{\mathbf M}_S=\left(\bigotimes\limits_{v=1}^{l-1}{\mathbf I}_{I_{q_v}}\right)
	\otimes {\mathbf M}_n\otimes
	\left(\bigotimes\limits_{v=l+1}^{N-K}{\mathbf I}_{I_{q_v}}\right)
	\label{eq:somepermutation}
		\end{equation}
	or  ${\mathbf M}_S =\mathbf I_K\otimes \mathbf M_n\otimes \mathbf I_L$,
	where $K=\prod\limits_{v=1}^{l-1}I_{q_v}$ and  $L=\prod\limits_{v=l+1}^{N-k}I_{q_v}$. 
	\end{lemma}
\begin{proof}
	Let $\delta(i,j)$ denote the Kronecker delta symbol, i.e., $\delta(i,j)=1$ for $i=j$ and $\delta(i,j)=0$ for $i\ne j$.
	One can easily verify that  
	\begin{equation}
	\begin{split}
	&\left(\mathbf M_S\right)_{
				ind_{i_{\tilde{q}_1},\dots,i_{\tilde{q}_{N-k}}}^{I_{q_1}\times\dots\times I_{q_{N-k}}}
			,\colvec{q}{N-k}} =\\
	& \delta(i_{\tilde{q}_1},i_{q_1})\cdots\delta(i_{\tilde{q}_{l-1}},i_{q_{l-1}})\cdot
	(\mathbf M_n)_{{\tilde{q}_l}, n}\cdot\delta(i_{\tilde{q}_{l+1}},i_{q_{l+1}})\cdots\delta(i_{\tilde{q}_{N-k}},i_{q_{N-k}})=\\
	&\begin{cases}
(\mathbf M_n)_{{\tilde{q}_l}, n},& \text{ if }
 i_{\tilde{q}_1}=i_{q_1},\dots, 
 i_{\tilde{q}_{l-1}}=i_{q_{l-1}},
 i_{\tilde{q}_{l+1}}=i_{q_{l+1}},\dots,
 i_{\tilde{q}_{N-k}}=i_{q_{N-k}}\\
 0,&\text{otherwise}.
	\end{cases}
	\end{split}
	\label{eq:matrixMs}
	\end{equation}   
	Hence
	\begin{gather*}
	\left(\unfset{\mathbf B}{S}\right)_{\colvec{p}{k},\colvec{q}{N-k}}\stackrel{\cref{eq:defmatrreprS}}{=}\left(\mathcal B\right)_{i_1,\dots,i_N}\stackrel{\cref{eq:defmatrreprS}}{=}
	\left(\unf{B}{n}\right)_{ind_{i_1,\dots,i_{n-1},i_{n+1},\dots,i_N}^{I_1\times\dots\times I_{n-1}\times I_{n+1}\dots\times I_N},i_n}=\\
	\left(\unf{A}{n}\mathbf M_n\right)_{ind_{i_1,\dots,i_{n-1},i_{n+1},\dots,i_N}^{I_1\times\dots\times I_{n-1}\times I_{n+1}\dots\times I_N},i_n}=\sum\limits_{\tilde{q}_l=1}^n (\mathcal A)_{i_1,\dots,i_{n-1},i_{\tilde{q}_l},i_{n+1},\dots,i_N}(\mathbf M_n)_{{\tilde{q}_l}, n}\stackrel{\cref{eq:defmatrreprS}}{=}\\
	\sum\limits_{\tilde{q}_l=1}^n \left(\unfset{\mathbf A}{S}\right)_{
	\colvec{p}{k}, 
ind_{i_{q_1},\dots,i_{q_{l-1}},i_{\tilde{q}_l},i_{q_{l+1}}\dots,i_{q_{N-k}}}^{I_{q_1}\times\dots\times I_{q_{N-k}}}
}
	(\mathbf M_n)_{{\tilde{q}_l}, n}\stackrel{\cref{eq:matrixMs}}{=}\\
	\sum\limits_{\tilde{q}_1=1}^{q_1}\cdots
	\sum\limits_{\tilde{q}_{N-k}=1}^{q_{N-k}}
		 \left(\unfset{\mathbf A}{S}\right)_{
		\colvec{p}{k}, 
		ind_{i_{\tilde{q}_1},\dots,i_{\tilde{q}_{N-k}}}^{I_{q_1}\times\dots\times I_{q_{N-k}}}
			}
	\left(\mathbf M_S\right)_{
		ind_{i_{\tilde{q}_1},\dots,i_{\tilde{q}_{N-k}}}^{I_{q_1}\times\dots\times I_{q_{N-k}}}
				,\colvec{q}{N-k}}=\\
	\left(\unfset{\mathbf A}{S}{\mathbf M}_S\right)_{\colvec{p}{k},\colvec{q}{N-k}}.
		\end{gather*}
\end{proof}
The following  corollary follows  from \cref{Thm:HOcase} and states that $2^N-2-N$ conditions in \cref{eq:manyinclusions} are redundant.
\begin{corollary}\label{corr:redandancy}
Let $N\geq 4$,  and let $\mathcal A, \mathcal B\in\mathbb F^{I_1\times \dots\times I_N}$. If
\cref{eq:AincludesBnew} holds, then  \cref{eq:manyinclusions} also holds.	
\end{corollary}
\section{Illustration: classification of linear mixtures of signals}\label{sec4:illustration} A basic problem in signal processing is to assess whether two observed signals involve the same underlying signal ``components''. Typically, the component signals manifest themselves with a different amplitude in the observed signals. If moreover the component signals are by themselves unknown, which is  the case in many applications, the problem can be very challenging. As a preview, in
\cref{fig:3} it may  a priori not be obvious to establish which displayed signals are generated by the same components up to scaling.  

One of the possible applications is in underdetermined 
Blind Source Separation (BSS). In BSS, the task is to recover  sources  from a set of their linear mixtures \cite{ComoJ10}.
Often, sources are sparsely combined in the observed mixed signals \cite{Donoho2001}, i.e.,
the number of sources is large but each mixture contains a small number of sources.
This means that the mixing matrix is sparse and has many more columns than rows. BSS problems that involve a wide mixing matrix
are called underdetermined and are generally much harder to solve than overdetermined BSS problems (involving a mixing matrix that is square or tall).
As a  preprocessing  step   one can first try to solve the following multi-label classification problem:  mixture $i$ belongs to the same class as mixture $j$ if  mixture $i$ is generated by (some of) the  sources  that appear in mixture $j$. In this way the initial underdetermined BSS problem with many sources can be  replaced by a set of smaller overdetermined BSS problems.

In this section we explain how  \cref{thm3rdorder} can be used to solve the  multi-label classification problem.
Our derivation is valid under the assumption that the sources can   simultaneously be mapped (i.e., ``tensorized'') into low ML rank tensors and that
the mapping, so called {\em tensorization}, is linear.
Such mappings are known \cite{OttoThesis,Martijn2017,Khoromskii2018} for sources that can be modeled as  exponential polynomials ({\em Hankelization}), rational functions ({\em L\"ownerization}), and  periodic signals ({\em Segmentation}), among others. To demonstrate the approach we confine ourselves to  exponential polynomials.

To solve the multi-label classification problem, we do not use more prior knowledge about the sources than that they can be (approximately) modeled as exponential polynomials (with a mild bound on the value $L_s$ in \cref{eq:Ls} that will be introduced in the next subsection). 

\subsection{Exponential polynomials and Hankelization mapping}\label{subsec:Hankelization}
A univariate exponential polynomial is a function of the form
\begin{equation}
s(t)=\sum\limits_{f=1}^Fp_f(t)a_f^t,
\label{eq:exp_poly}
\end{equation}
where $p_1,\dots,p_F$ are non-zero polynomials in one variable and $a_1,\dots,a_F\in\mathbb  C\setminus\{0\}$.
Let $T_s$ denote the sampling time and  let $N$ be the number of sampling points.
It can be  shown \cite{LievenLrLr1,OttoThesis} that for any positive integers 
$I_1,I_2,I_3$ that sum up to $N+2$ and are greater than or equal to $L_{\mathbf s}$, the vector $\mathbf s=[s_1\ \dots\ s_N]^T:=[s(0)\ \dots\ s((N-1)T_s)]^T$ can be  mapped to an $I_1\times I_2\times I_3$ ML rank-$(L_{\mathbf s},L_{\mathbf s},L_{\mathbf s})$ tensor $\mathcal S$, where the value 
\begin{equation}
L_{\mathbf s}:=F+\sum_{f=1}^F\deg p_f
\label{eq:Ls}
\end{equation}
does not depend on $I_1$, $I_2$, $I_3$. The mapping $H:\mathbf s\mapsto \mathcal S$, $H:\mathbb C^N\to \mathbb C^{I_1\times I_2\times I_3}$   is given by \cite{LievenLrLr1,OttoThesis}
$$
(\mathcal S)_{i_1i_2i_3}=s_{i_1+i_2+i_3-2}=s((i_1+i_2+i_3-3)T_s),
$$
where $1\leq i_1\leq I_1$, $1\leq i_2\leq I_2$, $1\leq i_3\leq I_3$.
Since $(\mathcal S)_{i_1i_2i_3}$ depends only on $i_1+i_2+i_3$, the mapping $H$ was called  ``Hankelization'' in \cite{OttoThesis}.
It is worth noting that if $I_1=I_2=I_3$, then $\mathcal S$ is a fully symmetric  tensor, implying
that $\unfone{S}=\unftwo{S}=\unfthree{S}$.

 It is clear that $H$ is a linear mapping, so if 
  $\mathbf y=[y_1\ \dots\ y_N]^T:=[y(0)\ \dots\ y((N-1)T_s)]^T$
 is a linear mixture of  sampled sources of the form \cref{eq:exp_poly}
\begin{equation}
y(t) = \mixentry_1s_1(t)+\dots \mixentry_Rs_r(t),\qquad t=0,T_s,\dots,(N-1)T_s
\label{eq:mixture}
\end{equation}
and $\min(I_1,I_2,I_3)\geq \max L_{{\mathbf s}_r}$, then, by \cref{eq:Ls},
\begin{equation}
\mathcal Y:=H(\mathbf y)=\mixentry_1H(\mathbf s_1)+\dots+\mixentry_RH(\mathbf s_R)=\mixentry_1\mathcal S_1+\dots+\mixentry_R\mathcal S_R
\label{eq:hankelizationofasum}
\end{equation}
is a decomposition of  $\mathcal Y$ into a sum of ML rank-$(L_{\mathbf s_r},L_{\mathbf s_r},L_{\mathbf s_r})$ terms.

\subsection{Example}\label{subsec:exampleclass}
\begin{figure}[htbp]
		\centering
		\label{fig:4}
		\includegraphics[trim=50 98 45 70, clip, width=0.7\textwidth]{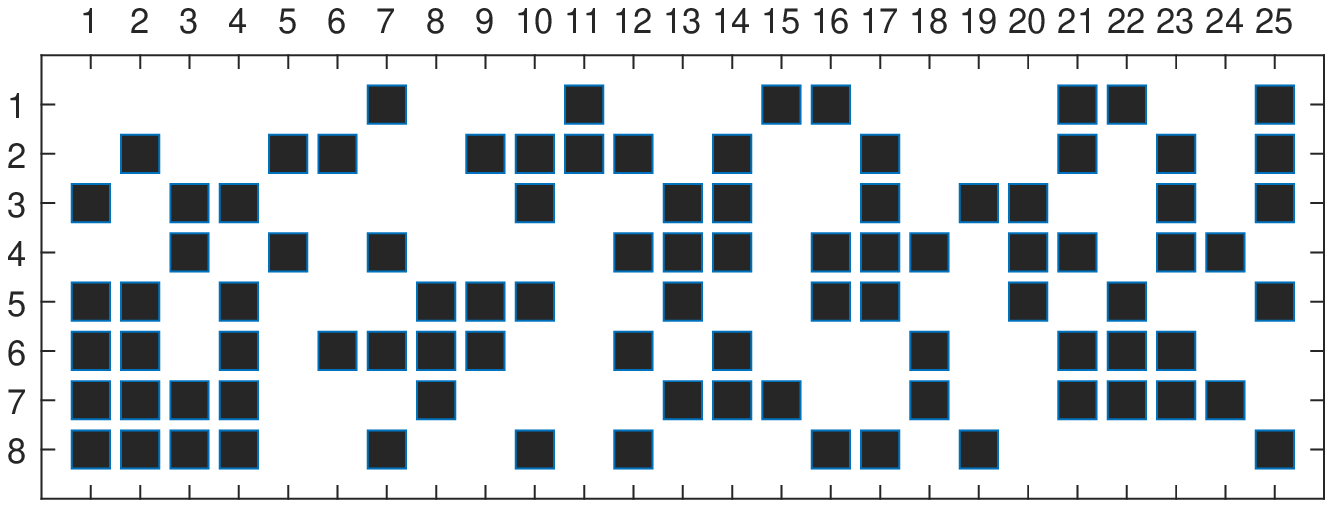}
		\caption{The nonzero pattern of the matrix $\mixmatrix$.
				}
		\label{fig:matrix_M}
	\end{figure}

\begin{figure}[htbp]
		\centering
		\label{fig:1}
		\includegraphics[scale=0.35]{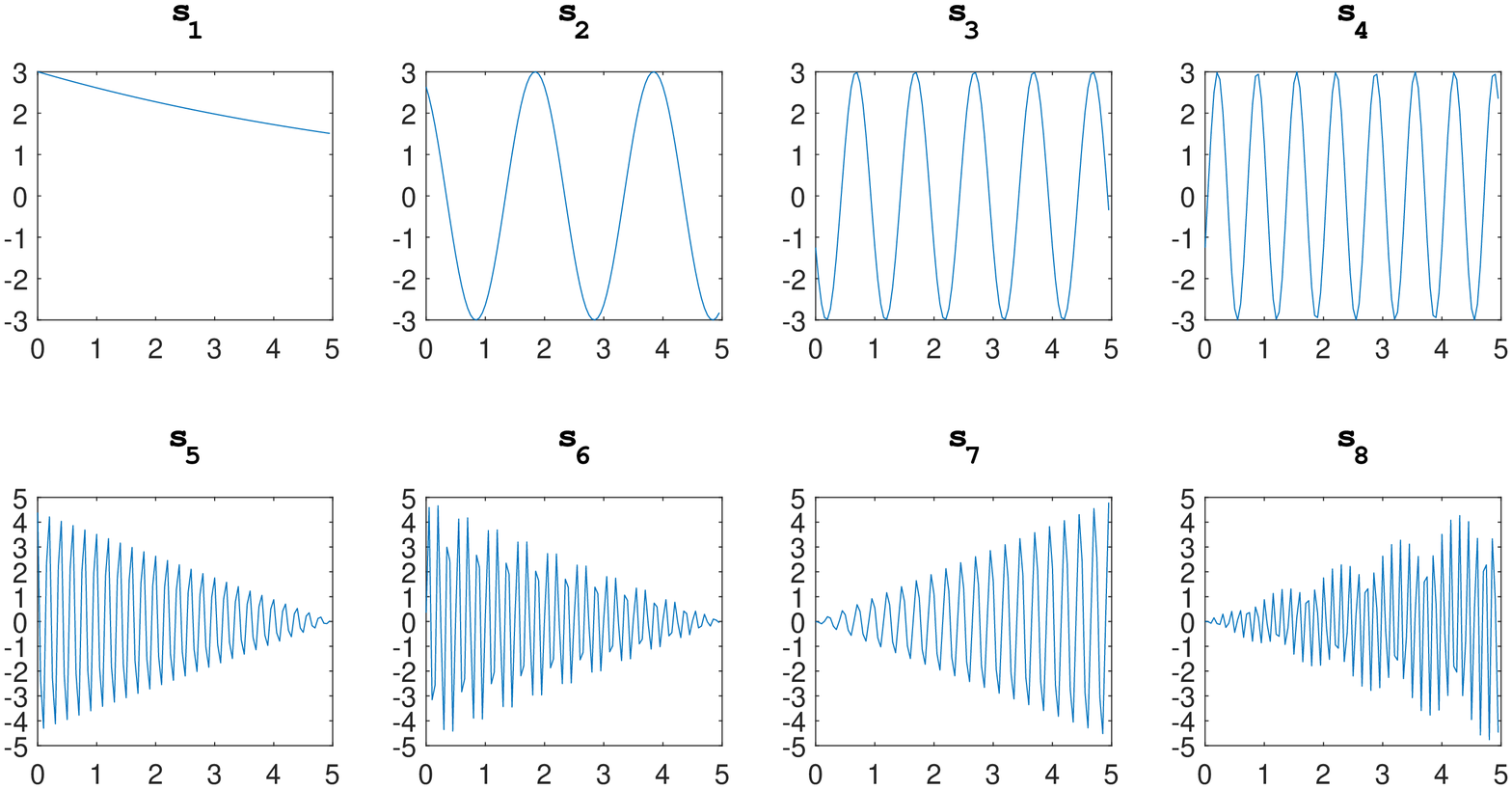}
		\caption{Sampled source signals ${\mathbf s}_1,\dots,{\mathbf s}_8$.
				}
		\label{fig:source_signals_time}
	\end{figure}

\begin{figure}[htbp]
		\centering
		\label{fig:3}
		\includegraphics[scale=0.35]{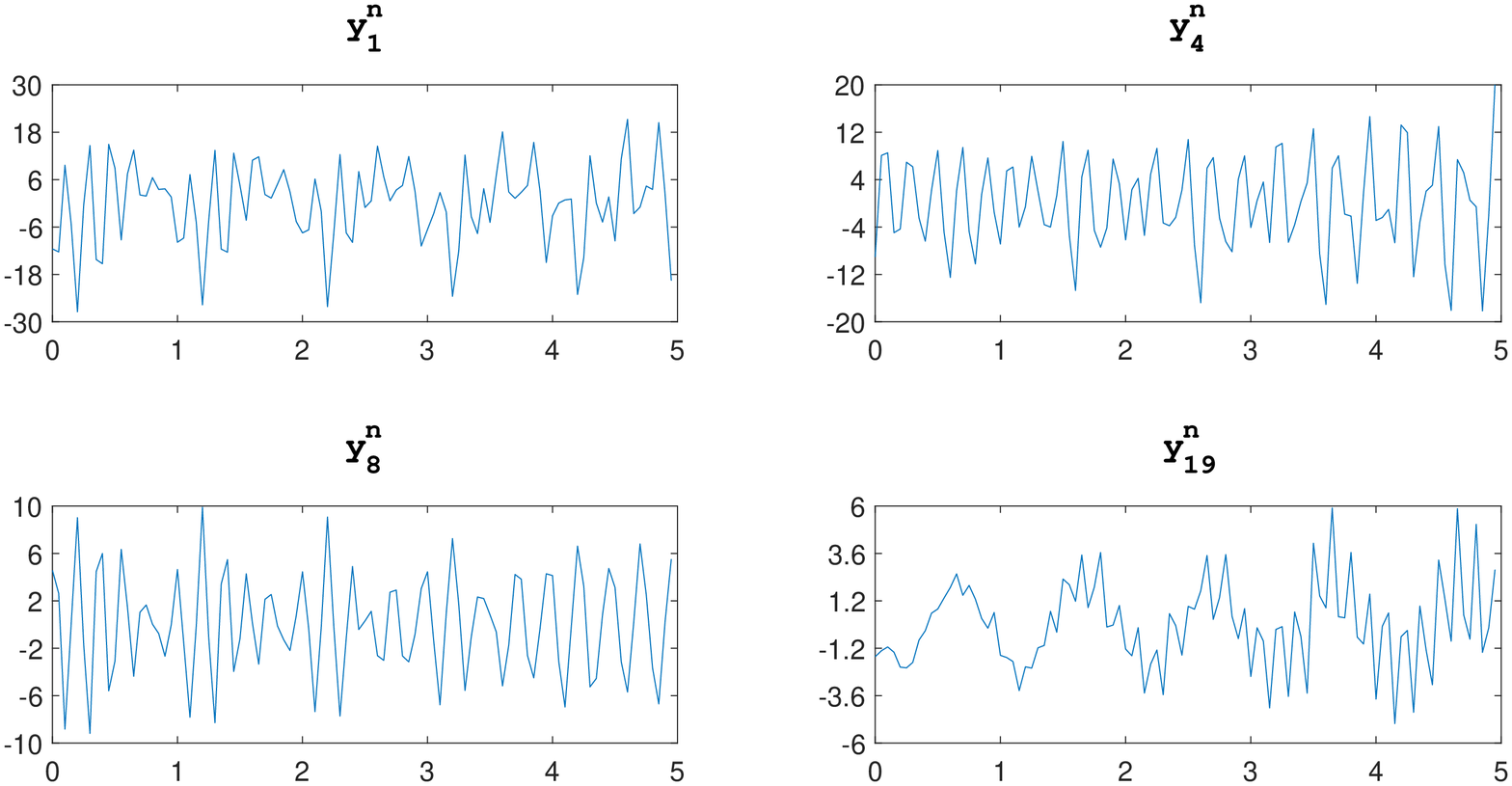}
		\caption{Noisy sampled mixtures
			${\mathbf y}_1^n$, ${\mathbf y}_4^n$, ${\mathbf y}_8^n$, ${\mathbf y}_{19}^n$ $($see \cref{eq:mixture1,eq:mixture4,eq:mixture8,eq:mixture19}$)$}.
				\label{fig:mixtures_8_19_1_14_time}
	\end{figure}

\begin{figure}[tbhp]
    	\centering
    	\subfloat[]{\label{fig:a}\includegraphics[trim=50 80 45 70, clip, width=0.5\textwidth]{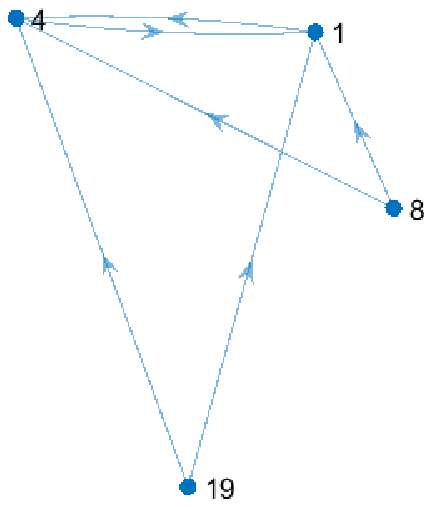}}
    	\subfloat[]{\label{fig:b}\includegraphics[trim=50 80 45 70, clip, width=0.5\textwidth]{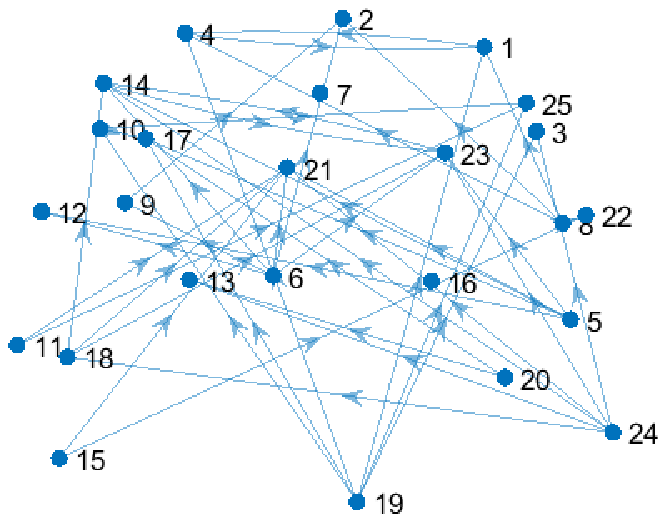}}
    	\caption{The subgraph corresponding to the mixtures 	$y_1(t)$, $y_4(t)$, $y_8(t)$, $y_{19}(t)$ $($left$)$ and the
    	graph corresponding to all mixtures $y_1(t),\dots,y_{25}(t)$ $($right$)$; a directed edge from  vertex $i$  to  vertex $j$ indicates that mixture $y_i(t)$ is generated by sources that also appear in  mixture $y_j(t)$.}
    	\label{fig:digraph}
    \end{figure}

	We generate $25$    mixtures 
	\begin{equation}
	y_j(t) = \mixentry_{1j}s_1(t)+\dots \mixentry_{8j}s_8(t),\qquad j=1,\dots,25
	\label{eq:sparse_mixture}
	\end{equation}
	of $8$ exponential polynomials
	  \begin{align*}
	&s_1(t)=3\cdot 2^{-\frac{t}{5}},&
	&\hspace{-1mm}s_2(t)=3\cos(\pi t+\frac{1}{2}),&
	&s_3(t)=3\cos(2\pi t+2),\\
	&s_4(t)=3\cos(3\pi t-2),&
	&\hspace{-1mm}s_5(t)=(5-t)\cos(10\pi t+\frac{1}{2}),&
	&s_6(t)=(5-t)\cos(12\pi t-\frac{3}{2}),
	\\
	&s_7(t)=t\cos(8\pi t+1),&
	&\hspace{-1mm}s_8(t)=t\cos(14\pi t-\frac{1}{2}).
	\end{align*}
	The coefficients $\mixentry_{ij}$ are generated randomly\footnote{The numerical experiments in the example were performed in MATLAB R2018b. To make  the results reproducible,  the random number generator was initialized	using the built-in function \texttt{rng('default')} (the Mersenne Twister with seed $0$).} so that for each $j=1,\dots,25$
	at least three and at most six of  $\mixentry_{1j},\dots,\mixentry_{8j}$ are zero. The nonzero
	coefficients $\mixentry_{ij}$   are randomly chosen from $[-2.5,-0.5]\cup [0.5, 2.5]$. 	We thus  obtain that 
	\begin{equation}
	[y_1(t)\ \dots  y_{25}(t)] = [ s_1(t)\ \dots s_{8}(t)]\mixmatrix,
	\label{Y=SMexact}
	\end{equation}
	where $\mixmatrix=(\mixentry_{ij})$ is an $8\times 25$ sparse matrix. The nonzero pattern of $\mixmatrix$ is shown in	\cref{fig:matrix_M}.
By way of example, the mixtures $y_1(t)$, $y_4(t)$, $y_8(t)$, and $y_{19}(t)$ were generated as
\begin{align}
&y_1(t) =&   2.22s_3(t) & &-1.95s_5(t)&	&-2.38s_6(t)& 	&-2.39s_7(t)\phantom{,}& &+1.77s_8(t),&\label{eq:mixture1}\\ 
&y_4(t) =& -0.55s_3(t)& &-2.07s_5(t)&	&+0.50s_6(t)& 	&+2.41s_7(t)\phantom{,}& &-1.90s_8(t),&\label{eq:mixture4}\\
&y_8(t) =&            & &+1.16s_5(t)&	&+0.94s_6(t)& 	&+1.35s_7(t),& & &\label{eq:mixture8}\\
&y_{19}(t) =&   0.69s_3(t)& &&	&& 	&& &-0.68s_8(t).&\label{eq:mixture19}
\end{align}
	We consider a noisy sampled (with $T_s=0.05$ and $N=100$)  version of \cref{Y=SMexact}:
	\begin{equation}
 [{\mathbf y}_1^n\ \dots\ {\mathbf y}_{25}^n]:=	[\mathbf y_1\ \dots\ \mathbf y_{25}]+\sigma\mathbf N = [\mathbf s_1\ \dots \mathbf s_{8}]\mixmatrix+\sigma\mathbf N,
	\label{eq:noisyand sampled}
	\end{equation}
	 in which the entries of the $100\times 25$ matrix $\mathbf N$ are  independently drawn from the standard normal distribution $N(0,1)$ and $\sigma=0.1\frac{\|[\mathbf y_1\ \dots\ \mathbf y_{25}]\|_F}{\|\mathbf N\|_F}$, where $\|\cdot\|_F$ denotes the Frobenius norm.\footnote{Note that the matrix pencil based algorithm in \cite{LievenLrLr1}
			can be used to estimate the matrix $\mixmatrix$ and the sources $\mathbf s_1,\dots,\mathbf s_8$ only for  much smaller values of $\sigma$.}
	The sampled sources  $\mathbf s_1,\dots,\mathbf s_8$ and  noisy sampled mixtures
	${\mathbf y}_1^n$, ${\mathbf y}_4^n$,
	 ${\mathbf y}_8^n$, ${\mathbf y}_{19}^n$ are shown in \cref{fig:source_signals_time,fig:mixtures_8_19_1_14_time}, respectively.
    We now  use \cref{thm3rdorder} to verify whether  the pair of mixtures $(y_i(t),y_j(t))$ is generated by the same subset of sources, $1\leq i,j\leq 25$\footnote{Note that    		
    		in contrast to BSS, we do not work with the full matrix $\mixmatrix$ but only with pairs of its columns. The number of mixtures $25$ is just chosen to illustrate the approach for a large number  of pairs (namely, $25^2-25=600$), making the results very convincing.}, $i\ne j$.   For visualization purposes it is convenient to associate the mixtures $y_1(t),\dots,y_{25}(t)$ with vertices of a  directed graph: a  directed edge from  vertex $i$  to  vertex $j$ indicates that  $y_i(t)$ is generated by sources that also appear in   $y_j(t)$.
    For instance, the subgraph corresponding to the mixtures $y_1(t)$, $y_4(t)$, $y_8(t)$, and $y_{19}(t)$ is shown in  \cref{fig:a}.
    In this example we show how  to recover the overall graph in  \cref{fig:b} based only on  the (observed) vectors ${\mathbf y}_1^n,\dots,{\mathbf y}_{25}^n$ and   without estimating ${\mathbf s}_1,\dots,{\mathbf s}_8$.
Let $H:\mathbb C^{100}\to \mathbb C^{34\times 34\times 34}$ denote the Hankelization mapping.
 By \cref{eq:hankelizationofasum,eq:noisyand sampled}, we have that
\begin{equation*}
{\mathcal Y}_i^n:=H({\mathbf y}_i^n)=\mixentry_{i1}H(\mathbf s_1)+\dots+\mixentry_{i8}H(\mathbf s_8)+\sigma H(\mathbf n_i)=\mixentry_{i1}\mathcal S_1+\dots+\mixentry_{18}\mathcal S_R+\sigma H(\mathbf n_i)
\end{equation*}
is an approximate decomposition of  ${\mathcal Y}_i^n$ into a sum of ML rank-$(L_{\mathbf s_r},L_{\mathbf s_r},L_{\mathbf s_r})$ terms.    
One can easily verify that the exact values of  $L_{\mathbf s_1},\dots,L_{\mathbf s_8}$ are 
$1,2,2,2,4,4,4,4$, respectively. For instance, since 
$$
s_8(t) = t\frac{1}{2}(e^{(14\pi t-0.5)i}+e^{-(14\pi t-0.5)i})=(\frac{1}{2}e^{-0.5i}t^1)e^{14\pi t}+(\frac{1}{2}e^{0.5i}t^1)e^{-14\pi t},
$$
we get, by \cref{eq:Ls}, that $L_{\mathbf s_8}=2+(1+1)=4$.
On the other hand, it can be verified that although the tensors $\mathcal S_5,\dots,\mathcal S_8$ are ML rank-$(4,4,4)$, they can be approximated by ML rank-$(2,2,2)$ tensors
with a relative error less than $0.061$, which is below the noise level.

To verify whether $y_i(t)$ is generated by sources that appear in $y_j(t)$,  it is sufficient to show that
${\mathcal Y}_i^n$ is generated by ML rank terms that appear in ${\mathcal Y}_j^n$, which, by 
 \cref{thm3rdorder}, is reduced to  verifying  that the column space of $\unfgenone{\mathbf Y_i^n}$
 is contained in the column space of $\unfgenone{\mathbf Y_j^n}$. To compare the column spaces we proceeded as follows. 
 For each $i=1,\dots,25$ we computed the first $r_i$  singular vectors $\mathbf u_{i1},\dots,\mathbf u_{ir_i} $  of  $\unfgenone{\mathbf Y_i^n}$, where
 the rank $r_i$ of $\unfgenone{\mathbf Y_i^n}$ was estimated   as the largest index $k$ such that
 the ratio of the $k$th and the $(k+1)$st singular value of $\unfgenone{\mathbf Y_i^n}$ is greater than a certain threshold $\tau$. (We chose  $\tau=2.3$.)
 Then for all  $i,j=1,\dots,25$ and $i\ne j$ we concluded that
 the column space of $\unfgenone{\mathbf Y_i^n}$  is contained in the column space of $\unfgenone{\mathbf Y_j^n}$  if
 the $(r_j+1)$st singular value of the matrix $[\mathbf u_{i1}\ \dots\ \mathbf u_{ir_i}\  \mathbf u_{j1}\ \dots\ \mathbf u_{jr_j}]$
 was less than $0.1$ (and in  this case we plotted the directed edge from  vertex $i$ to  vertex $j$). The resulting directed graph is shown in
 \cref{fig:b}. The same graph can be obtained directly from the nonzero pattern of the matrix $\mixmatrix$ which means that all $43$ (out of the possible $600$) edges of the graph were detected correctly and no superfluous edges were added.
 \section{Conclusion}
  An obvious requirement for a tensor $\mathcal B$ to be the sum of (possibly scaled) terms from the decomposition of a tensor $\mathcal A$, is that its
 column (row, fiber, \dots) space is a subspace of the corresponding space of tensor $\mathcal A$. Formally, this means that
 $\operatorname{row}(\unfset{\mathbf B}{n}) \subseteq \operatorname{row}(\unfset{\mathbf A}{n})$
 should hold for all $n\in\{1,\dots, N\}$. However, this is only a \textit{necessary} condition. Switching to the column spaces, we have shown in this paper that 
 \begin{equation}
 \operatorname{col}(\unfset{\mathbf B}{n}) \subseteq \operatorname{col}(\unfset{\mathbf A}{n}),\qquad n\in\{1,\dots, N\}
 \label{eq:lastagain}
 \end{equation}
  is a \textit{sufficient} condition for $\mathcal B$ to be generated by (possibly scaled) terms from the decomposition of $\mathcal A$.
  The number or terms and their ``type'' (namely, their ML rank)  follow from the analysis as well. As the derivation relies only on linear algebra, it bypasses the typical difficulties in the computation of CPD and BTD, such as NP-hardness and possible ill-conditioning. We believe that this paper introduces a new tool that will prove important for tensor-based pattern recognition and machine learning, in a similar way as (explicit)  tensor decompositions have proven to be fundamental tools for data analysis. We have illustrated the practical use of the new tool in a new clustering-based scheme for sparse underdetermined BSS.
  
   An interesting topic of further study would be to investigate  partially shared structure of $\mathcal A$ and $\mathcal B$, in the sense that $\mathcal A$ and $\mathcal B$ share some but not all terms.
We will also derive more detailed information from the actual principal angles  and associated directions between the subspaces obtained from $\mathcal A$ and $\mathcal B$. Another topic of further study is the generalization to ``flower'', ``butterfly'' and related decompositions \cite{Martijn20172,Martijn2017,Sorber}.
\bibliographystyle{siamplain}
\bibliography{tens_sim}
\end{document}